\documentclass{article}
\usepackage[latin9]{inputenc}
\usepackage{mathtools}
\usepackage{amsmath}
\usepackage{amssymb}

\makeatletter
\usepackage{amsfonts}
\usepackage{graphicx}
\providecommand{\U}[1]{\protect\rule{.1in}{.1in}}
\newtheorem{theorem}{Theorem}\newtheorem{definition}[theorem]{Definition}\newtheorem{lemma}[theorem]{Lemma}\newtheorem{proposition}[theorem]{Proposition}\newtheorem{remark}[theorem]{Remark}
\makeatother

\begin{document}

\title{On the relationship between the \emph{energy shaping} and the \emph{Lyapunov
constraint based} methods}

\author{Sergio Grillo,$^{a,b}$ Leandro Salomone$^{b,c}$\ \ \&\ \ \ Marcela
Zuccalli$^{c}$\\
$^{a}$\textit{Centro Atómico Bariloche and Instituto Balseiro}, \textit{8400
S.C. de Bariloche, Argentina}\\
$^{b}$CONICET\textit{, Argentina}\\
$^{c}$\textit{Departamento de Matemática,} \textit{Facultad de Ciencias
Exactas, U.N.L.P., Argentina}}
\maketitle
\begin{abstract}
In this paper, we make a review of the controlled Hamiltonians (CH)
method and its related matching conditions, focusing on an improved
version recently developed by D.E. Chang. Also, we review the general
ideas around the Lyapunov constraint based (LCB) method, whose related
partial differential equations (PDEs) were originally studied for
underactuated systems with only one actuator, and then we study its
PDEs for an arbitrary number of actuators. We analyze and compare
these methods within the framework of Differential Geometry, and from
a purely theoretical point of view. We show, in the context of underactuated
systems defined by simple Hamiltonian functions, that the LCB method
and the Chang's version of the CH method are equivalent stabilization
methods (i.e. they give rise to the same set of control laws). In
other words, we show that the Chang's improvement of the energy shaping
method is precisely the LCB method. As a by-product, coordinate-free
and connection-free expressions of Chang's matching conditions are
obtained.
\end{abstract}

\section{Introduction}

Under the name of \emph{energy shaping method}, several methods or
procedures for achieving (asymptotic) stabilization of nonlinear underactuated
Lagrangian and Hamiltonian systems are included: \emph{potential shaping},
\emph{kinetic shaping}, \emph{total energy shaping}, \emph{energy
plus force shaping}, \emph{IDA-PBC}, etc. See for instance \cite{auck,ks,cl,es1,es2,kri,idapbc,shaft2,woosley},
and \cite{chang4,romero} for more recent works. They are based on
the idea of \emph{feedback equivalence} (see Ref. \cite{CBLMW02}),
and their purpose is to construct, for a given underactuated mechanical
system, a control law and a Lyapunov function for the resulting closed-loop
system. To do that, a set of partial differential equations (PDEs),
known as \emph{matching conditions}, must be solved. Such PDEs have
among their unknowns the aforementioned Lyapunov function.

All of these methods can be seen as particular versions of the so-called
\emph{controlled Lagrangians (CL) method} or the \emph{controlled
Hamiltonians (CH) method}, which in turn are equivalent stabilization
methods, in a sense that has been carefully explained in Ref. \cite{CBLMW02}.

The origin of the energy shaping method can be placed 35 years ago
\cite{ari,bro,shaft,will}, while the method in its more general form
is around 15 years old \cite{CBLMW02}. More recently, 6 years ago,
an alternative stabilization method for nonlinear underactuated mechanical
systems has been presented: the \emph{Lyapunov constraint based (LCB)
method}. It appeared for the first time in \cite{hocs}, it was further
developed in \cite{gym}, and it was extended to systems with impulsive
effects in Ref. \cite{chaalal}. The method is based on the idea of
controlling actuated mechanical systems by imposing kinematic constraints
(see \cite{cg,g,marle,marle2,perez1,perez2,shiriaev}). It serves
the same purpose as the energy shaping method (to construct a control
law and a Lyapunov function for the resulting closed-loop system)
and, in order to accomplish it, a set of PDEs must be solved too.
It is worth mentioning that the LCB method has been originally developed
for underactuated systems with only one actuator.

One of the aims of this paper is to extend the study of the LCB method
to an arbitrary number of actuators, and to show that this method
contains every version of the energy shaping method (and actually,
every method that serves the same purpose), in the sense that the
set of control laws that can be constructed with the energy shaping
method is contained in the corresponding set of the LCB method (extending
a result already presented in \cite{gym}).

Almost simultaneously with the appearance of \cite{gym}, an improvement
of the energy (plus force) shaping method, for underactuated systems
defined by \emph{simple} Lagrangian or Hamiltonian functions, was
presented by Chang in \cite{chang1,chang2,chang3}. It consists in
an important simplification of the matching conditions. The main goal
of the present paper is to show that such matching conditions are
exactly the PDEs related to the LCB method, at least in the context
of simple Hamiltonian functions. Moreover, we show in the same context
that the Chang's version of the energy shaping method is equivalent
to the LCB method, i.e. both methods give rise exactly, for a given
underactuated system, to the same set of control laws. In other words,
we show that the Chang's improvement of the energy shaping method
is precisely the LCB method. Such a result is quite surprising for
us, because the involved methods are based on very different ideas:
``feedback equivalence'' and ``controlling by the imposition of
kinematic constraints.'' 

We can say that this article is similar in spirit to Ref. \cite{CBLMW02},
where the equivalence between the CL and the CH methods was established.
In particular, as in that paper, a substantial portion of the work
is dedicated to describe, in a very precise way and by using the same
language, the methods that we want to compare.

\bigskip{}

The paper is organized as follows. In Section 2, we present some basic
facts about affine connections on general linear bundles, which will
be used along all of the paper to write down coordinate-free expressions
of the PDEs that we want to study. In Section 3, we make a review
of the energy shaping method in a Hamiltonian language, i.e. the controlled
Hamiltonians (CH) or IDA-PBC method. We begin with a rather general
version of the method, then we progressively consider particular situations,
and finally we present Chang's version of the method (see for instance
\cite{chang4}), with its related matching conditions. In Section
4, we recall the idea of controlling mechanical systems by the imposition
of kinematic constraints. In particular, we review the idea of achieving
(asymptotic) stability by means of the so-called \emph{Lyapunov constraint},
which gives rise to the LCB method and its related set of PDEs. We
show in the last section of the paper that such PDEs are exactly the
matching conditions obtained by Chang \cite{chang4}, at least when
underactuated systems defined by simple Hamiltonian functions are
considered. Finally, we show the equivalence of the LCB method and
the Chang's version of the CH method.

\bigskip{}

We assume that the reader is familiar with basic concepts of Differential
Geometry \cite{boot,kn,mrgm}, Hamiltonian systems in the context
of Geometric Mechanics \cite{am,ar,mr}, Control Theory in a geometric
language \cite{bloc,bullo}, and Lyapunov theorems for (asymptotic)
stability \cite{khalil}.

\bigskip{}

\textbf{Basic notation. }Along all of the paper, every manifold will
be a smooth finite dimensional manifold, typically denoted by $Q$.
By $\tau_{Q}:TQ\rightarrow Q$ and $\pi_{Q}:T^{\ast}Q\rightarrow Q$
we will denote the tangent and cotangent vector bundles, respectively.
As it is customary, we indicate by $\left\langle \cdot,\cdot\right\rangle $
the natural pairing between $T_{q}^{\ast}Q$ and $T_{q}Q$ at every
$q\in Q$, and by $\mathfrak{X}\left(Q\right)$ and $\Omega^{1}\left(Q\right)$
the sheaves of sections of $\tau_{Q}$ and $\pi_{Q}$, respectively.
Unless a confusion may arise, we shall omit the subindex $Q$ for
$\tau_{Q}$ and $\pi_{Q}$. For a vector field $Y:Q\rightarrow TQ$,
in order to indicate that its image is contained inside some subset
$W$ of $TQ$, we shall write, for simplicity, $Y\subset W$. Given
a second manifold $P$ and a smooth function $F:Q\rightarrow P$,
we denote by $F_{\ast}$ and $F^{\ast}$ the push-forward map and
its transpose, respectively. 

Consider a local chart $\left(U,\varphi\right)$ of $Q$, with $\varphi:U\rightarrow\mathbb{R}^{n}$.
Given $q\in U$, we write $\varphi\left(q\right)=\left(q^{1},...,q^{n}\right)=\mathbf{q}$.
For the induced local charts $\left(TU,\varphi_{\ast}\right)$ and
$\left(T^{\ast}U,\left(\varphi^{\ast}\right)^{-1}\right)$ on $TQ$
and $T^{\ast}Q$, respectively, we write
\begin{equation}
\begin{array}{l}
\varphi_{\ast}\left(v\right)=\left(q^{1},...,q^{n},\dot{q}^{1},...,\dot{q}^{n}\right)=\left(\mathbf{q},\mathbf{\dot{q}}\right),\\
\\
\left(\varphi^{\ast}\right)^{-1}\left(\alpha\right)=\left(q^{1},...,q^{n},p_{1},...,p_{n}\right)=\left(\mathbf{q},\mathbf{p}\right),
\end{array}\label{icor}
\end{equation}
or simply 
\begin{equation}
\varphi_{\ast,q}\left(v\right)=\mathbf{\dot{q}}\ \ \ \text{and\ \ \ }\left(\varphi_{q}^{\ast}\right)^{-1}\left(\alpha\right)=\mathbf{p},\label{icor1}
\end{equation}
for all $v\in TU$ and $\alpha\in T^{\ast}U$. On $TT^{\ast}Q$ we
shall consider the induced charts $\left(TT^{\ast}U,\left(\varphi^{\ast}\right)_{\ast}^{-1}\right)$,
and write 
\begin{equation}
\left(\varphi^{\ast}\right)_{\ast}^{-1}\left(V\right)=(\mathbf{q},\mathbf{p},\mathbf{\dot{q}},\mathbf{\dot{p}}),\label{icor2}
\end{equation}
for all $V\in TT^{\ast}U$.

\section{Some preliminary results}

In this section we shall recall some results on vector bundles and
affine connections that will enable us to write global expressions
of the equations we want to study later. Most of these results were
proved in Ref. \cite{g}. Nevertheless, for the sake of completeness,
we include some proofs here. Also, at the end of the section, we recall
some basic facts about Lyapunov functions.

\bigskip{}

Let us consider a vector bundle $\Pi:\mathcal{U}\rightarrow Q$ and
fix an affine connection $\nabla:\mathfrak{X}\left(Q\right)\times\Gamma\left(\mathcal{U}\right)\rightarrow\Gamma\left(\mathcal{U}\right)$.
Related to the latter we can define a diffeomorphism $\beta:T\mathcal{U}\rightarrow\mathcal{U}\oplus TQ\oplus\mathcal{U}$,
given as follows (see Ref. \cite{g}). For $V\in T\mathcal{U}$, consider
a curve $u:\left(-\varepsilon,\varepsilon\right)\rightarrow\mathcal{U}$
passing through $\tau_{\mathcal{U}}\left(V\right)$ and with velocity
$V$ at $s=0$, i.e. $u_{\ast}\left(\left.d/ds\right\vert _{0}\right)=V$.
Finally, define 
\[
\beta\left(V\right):=\tau_{\mathcal{U}}\left(V\right)\oplus\Pi_{\ast}\left(V\right)\oplus\frac{D}{Ds}u\left(0\right).
\]
Fixing $q\in Q$ and a vector $X\in\mathcal{U}_{q}$ (i.e. $\Pi\left(X\right)=q$),
we have the linear isomorphisms
\begin{equation}
\beta_{X}:T_{X}\mathcal{U}\rightarrow T_{q}Q\oplus\mathcal{U}_{q}\ \ \ \text{and}\ \ \ \beta_{X}^{-1}:T_{q}Q\oplus\mathcal{U}_{q}\rightarrow T_{X}\mathcal{U},\label{12aa}
\end{equation}
given by
\begin{equation}
\beta_{X}\left(V\right):=\Pi_{\ast}\left(V\right)\oplus\frac{D}{Ds}u\left(0\right)\ \ \ \ \ \text{and}\ \ \ \ \ \beta_{X}^{-1}\left(Y\oplus Z\right):=u_{\ast}\left(\left.d/ds\right\vert _{0}\right),\label{12a}
\end{equation}
respectively, where we take $u$ in the second equation to be a curve
in $\mathcal{U}$ such that 
\[
u(0)=X,\qquad\left(\Pi\circ u\right)_{\ast}(d/ds|_{0})=Y\qquad\text{and}\qquad\frac{D}{Ds}u(0)=Z.
\]
We have in addition their corresponding transpose maps
\begin{equation}
\beta_{X}^{\ast}:T_{q}^{\ast}Q\oplus\mathcal{U}_{q}^{\ast}\rightarrow T_{X}^{\ast}\mathcal{U}\ \ \ \ \text{and}\ \ \ \ \ \beta_{X}^{\ast-1}:T_{X}^{\ast}\mathcal{U}\rightarrow T_{q}^{\ast}Q\oplus\mathcal{U}_{q}^{\ast}.\label{12bb}
\end{equation}
In terms of $\beta$, the horizontal and vertical subbundles related
to $\nabla$ at a point $X\in\mathcal{U}_{q}$ are, respectively,
\[
\operatorname{Hor}_{X}:=\beta_{X}^{-1}\left(T_{q}Q\oplus0\right)\ \ \ \text{and}\ \ \ \operatorname{Ver}_{X}:=\ker\Pi_{\ast,X}=\beta_{X}^{-1}\left(0\oplus\mathcal{U}_{q}\right).
\]
It can be shown that the \emph{vertical lift isomorphism 
\[
\operatorname{vlift}_{X}:\mathcal{U}_{q}\rightarrow\ker\Pi_{\ast,X}:Z\longmapsto\left.\frac{d}{ds}\right\vert _{0}\left(X+s\,Z\right)
\]
}is related with $\beta_{X}^{-1}$ by the formula
\begin{equation}
\operatorname{vlift}_{X}\left(Z\right)=\beta_{X}^{-1}\left(0\oplus Z\right).\label{vb}
\end{equation}
(This is true for any connection $\nabla$.) Let us suppose that $\mathcal{U}=T^{\ast}Q$
and $\Pi=\pi_{Q}=\pi$. The related diffeomorphism
\begin{equation}
\beta:TT^{\ast}Q\rightarrow T^{\ast}Q\oplus TQ\oplus T^{\ast}Q\label{beta}
\end{equation}
is given by
\[
\beta(V)=\alpha\oplus\pi_{\ast}(V)\oplus\frac{D}{Ds}u(0),\ \ \forall\alpha\in T^{\ast}Q,\ \ V\in T_{\alpha}T^{\ast}Q,
\]
where $u:(-\varepsilon,\varepsilon)\rightarrow T^{\ast}Q$ is a curve
passing through $\alpha$ at $s=0$ with velocity $V$. In a local
chart $\left(U,\varphi\right)$ of $Q$, it is easy to show that
\begin{equation}
\beta(\mathbf{q},\mathbf{p},\mathbf{\dot{q}},\mathbf{\dot{p}})=\left(\mathbf{q},\mathbf{p}\right)\oplus\left(\mathbf{q},\mathbf{\dot{q}}\right)\oplus\left(\mathbf{q},\mathbf{\dot{p}}+\Gamma\left(\mathbf{q},\mathbf{p},\mathbf{\dot{q}}\right)\right)\label{bcpq}
\end{equation}
(omitting in the last expression the map $\varphi$, just for simplicity),
where $\Gamma\left(\mathbf{q},\mathbf{p},\mathbf{\dot{q}}\right)$
is given by the Christoffel symbols $\Gamma_{il}^{k}\left(\mathbf{q}\right)$
of $\nabla$ (in the coordinate frame) as 
\[
\Gamma_{i}\left(\mathbf{q},\mathbf{p},\mathbf{\dot{q}}\right)=\Gamma_{il}^{k}\left(\mathbf{q}\right)\,\ p_{k}\,\dot{q}^{l}.
\]
Sum over repeated indices convention is assumed from now on. On the
other hand, using the relationship between the vertical lift isomorphism
$\operatorname{vlift}_{\alpha}:T_{\pi(\alpha)}^{\ast}Q\rightarrow\operatorname{ker}\pi_{\ast,\alpha}$
and the linear isomorphism $\beta_{\alpha}:T_{\alpha}T^{\ast}Q\rightarrow T_{\pi\left(\alpha\right)}Q\oplus T_{\pi\left(\alpha\right)}^{\ast}Q$
{[}see Eqs. $\left(\ref{12aa}\right)$, $\left(\ref{12a}\right)$
and $\left(\ref{vb}\right)${]}, every vertical vector $Y_{\alpha}\in T_{\alpha}T^{\ast}Q$
may be identified with a unique covector $y_{\alpha}\in T_{\pi(\alpha)}^{\ast}Q$
in the following ways:
\begin{equation}
Y_{\alpha}=\operatorname{vlift}_{\alpha}\left(y_{\alpha}\right)=\beta_{\alpha}^{-1}(0\oplus y_{\alpha})=\beta^{-1}(\alpha\oplus0\oplus y_{\alpha}).\label{Yconbeta}
\end{equation}
As a consequence, every vertical vector field $Y:T^{\ast}Q\rightarrow TT^{\ast}Q$
is defined by the unique fiber preserving map $y:T^{\ast}Q\rightarrow T^{\ast}Q$
such that
\begin{equation}
Y\left(\alpha\right)=\operatorname{vlift}_{\alpha}\left(y\left(\alpha\right)\right)=\beta^{-1}(\alpha\oplus0\oplus y\left(\alpha\right)).\label{fieldconbeta}
\end{equation}

\begin{definition} Given a function $F:\mathcal{U}\rightarrow\mathbb{R}$,
the \textbf{fiber} and \textbf{base derivatives} of $F$ are defined
as the fiber-preserving maps $\mathbb{F}F:\mathcal{U}\rightarrow\mathcal{U}^{\ast}$
and $\mathbb{B}F:\mathcal{U}\rightarrow T^{\ast}Q$ given by 
\begin{equation}
\left\langle \mathbb{F}F\left(X\right),Z\right\rangle =\left.\frac{d}{ds}F\left(X+s\,Z\right)\right\vert _{s=0}\label{fiberder}
\end{equation}
and
\[
\left\langle \mathbb{B}F\left(X\right),Y\right\rangle =\left.\frac{d}{ds}F\left(u\left(s\right)\right)\right\vert _{s=0},
\]
respectively, where $u:\left(-\varepsilon,\varepsilon\right)\rightarrow\mathcal{U}$
is a (horizontal)\ curve such that
\begin{equation}
u\left(0\right)=X,\ \ \ \left(\Pi\circ u\right)_{\ast}\left(\left.d/ds\right\vert _{0}\right)=Y\ \ \ \text{and}\ \ \ \frac{D}{Ds}u\left(s\right)=0.\label{w2}
\end{equation}

\end{definition}

\begin{remark} Note that $\mathbb{F}F$ is independent of $\nabla$,
but $\mathbb{B}F$ is not. \end{remark}

Let us come back to the cotangent bundle of $Q$. Given a smooth function
$F:T^{\ast}Q\rightarrow\mathbb{R}$, the fiber and base derivatives
of $F$ are bundle morphisms $\mathbb{F}F:T^{\ast}Q\rightarrow TQ$
and $\mathbb{B}F:T^{\ast}Q\rightarrow T^{\ast}Q$, respectively, which
in local coordinates read
\begin{equation}
\left(\mathbb{F}F\left(\mathbf{q},\mathbf{p}\right)\right)^{i}=\frac{\partial F}{\partial p_{i}}\left(\mathbf{q},\mathbf{p}\right)\label{lfd}
\end{equation}
and
\begin{equation}
\left(\mathbb{B}F\left(\mathbf{q},\mathbf{p}\right)\right)_{i}=\frac{\partial F}{\partial q^{i}}\left(\mathbf{q},\mathbf{p}\right)+\Gamma_{il}^{k}\left(\mathbf{q}\right)\ \frac{\partial F}{\partial p_{l}}\left(\mathbf{q},\mathbf{p}\right)\,\,p_{k}.\label{lbd}
\end{equation}

\bigskip{}
Regarding \emph{basic functions} $F:\mathcal{U}\rightarrow\mathbb{R}$,
i.e. functions for which there exists $f:Q\rightarrow\mathbb{R}$
such that $F=f\circ\Pi$, we have the next result.

\begin{proposition} \label{bas}If $F:\mathcal{U}\rightarrow\mathbb{R}$
is basic, then
\begin{equation}
\mathbb{F}F=0\ \ \ \text{and\ \ \ }\mathbb{B}F=df\circ\Pi.\label{basic}
\end{equation}

\end{proposition}

\emph{Proof.} Given $X,Z\in\mathcal{U}_{q}$ for some $q\in Q$, i.e.
$\Pi\left(X\right)=\Pi\left(Z\right)=q$, we have that
\[
\left\langle \mathbb{F}F\left(X\right),Z\right\rangle =\left.\frac{d}{ds}F\left(X+s\,Z\right)\right\vert _{s=0}=\left.\frac{d}{ds}f\circ\Pi\left(X+s\,Z\right)\right\vert _{s=0}=\left.\frac{d}{ds}f\left(q\right)\right\vert _{s=0}=0.
\]
On the other hand, given in addition $Y\in T_{q}Q$ and a curve $u$
satisfying $\left(\ref{w2}\right)$,
\[
\left\langle \mathbb{B}F\left(X\right),Y\right\rangle =\left.\frac{d}{ds}F\left(u\left(s\right)\right)\right\vert _{s=0}=\left.\frac{d}{ds}f\left(\Pi\left(u\left(s\right)\right)\right)\right\vert _{s=0}=\left\langle df\left(\Pi\left(X\right)\right),Y\right\rangle ,
\]
as we wanted to show.\ \ \ $\square$

\bigskip{}

The isomorphisms $\beta_{X}^{\ast-1}$ {[}see Eq. $\left(\ref{12bb}\right)${]}
give rise to another diffeomorphism $\widetilde{\beta}:T^{\ast}\mathcal{U}\rightarrow\mathcal{U}\oplus T^{\ast}Q\oplus\mathcal{U}^{\ast},$
being $\widetilde{\beta}\left(\Sigma\right)=\beta_{X}^{\ast-1}\left(\Sigma\right)$
for all $X\in\mathcal{U}$ and $\Sigma\in T_{X}^{\ast}\mathcal{U}$.
For the cotangent bundle, we have a diffeomorphism
\begin{equation}
\widetilde{\beta}:T^{\ast}T^{\ast}Q\rightarrow T^{\ast}Q\oplus T^{\ast}Q\oplus TQ.\label{betam}
\end{equation}

\begin{proposition} Given $F:\mathcal{U}\rightarrow\mathbb{R}$ and
$X\in\mathcal{U}$,
\begin{equation}
\widetilde{\beta}\left(dF\left(X\right)\right)=X\oplus\mathbb{B}F\left(X\right)\oplus\mathbb{F}F\left(X\right).\label{bt}
\end{equation}

\end{proposition}

\emph{Proof.} We must show that, for all $q\in Q$, $Y\in T_{q}Q$
and $Z\in\mathcal{U}_{q}$, 
\[
\left\langle \beta_{X}^{\ast-1}\left(dF\left(X\right)\right),Y\oplus Z\right\rangle =\left\langle \mathbb{B}F\left(X\right),Y\right\rangle +\left\langle \mathbb{F}F\left(X\right),Z\right\rangle ,
\]
or equivalently, $\left\langle dF\left(X\right),V\right\rangle =\left\langle \mathbb{B}F\left(X\right),Y\right\rangle +\left\langle \mathbb{F}F\left(X\right),Z\right\rangle $,
for $V=\beta_{X}^{-1}\left(Y\oplus Z\right)$. Let $u_{1}:\left(-\varepsilon,\varepsilon\right)\rightarrow\mathcal{U}$
be a curve satisfying $\left(\ref{w2}\right)$ and $u_{2}:\left(-\varepsilon,\varepsilon\right)\rightarrow\mathcal{U}$
such that $u_{2}\left(s\right):=X+s\,Z$. Since 
\[
\beta_{X}\left(\left(u_{1}\right)_{\ast}\left(\left.d/ds\right\vert _{0}\right)\right)=Y\oplus0\;\;\;\text{and}\;\;\;\beta_{X}\left(\left(u_{2}\right)_{\ast}\left(\left.d/ds\right\vert _{0}\right)\right)=0\oplus Z,
\]
then $\left(u_{1}\right)_{\ast}\left(\left.d/ds\right\vert _{0}\right)+\left(u_{2}\right)_{\ast}\left(\left.d/ds\right\vert _{0}\right)=V$.
Consequently, 
\[
\left\langle dF\left(X\right),V\right\rangle =\left.\frac{d}{ds}F\left(u_{1}\left(s\right)\right)\right\vert _{s=0}+\left.\frac{d}{ds}F\left(X+s\,Z\right)\right\vert _{s=0},
\]
what ends our proof.\ \ \ $\square$

\bigskip{}

\begin{remark} Let us replace $\mathcal{U}$ by the Whitney sum of
$k$ copies of $\mathcal{U}$, which we shall denote $\mathcal{U}\times\cdot\cdot\cdot\times\mathcal{U}$,
and consider on such a vector bundle the affine connection naturally
induced by one fixed on $\mathcal{U}$. Then, given a function $F:\mathcal{U}\times\cdot\cdot\cdot\times\mathcal{U}\rightarrow\mathbb{R}$,
its fiber and base derivatives
\[
\mathbb{F}F:\mathcal{U}\times\cdot\cdot\cdot\times\mathcal{U}\rightarrow\mathcal{U}^{\ast}\times\cdot\cdot\cdot\times\mathcal{U}^{\ast}\ \ \ \text{and\ \ \ }\mathbb{B}F:\mathcal{U}\times\cdot\cdot\cdot\times\mathcal{U}\rightarrow T^{\ast}Q
\]
will be defined by the formulae 
\[
\left\langle \mathbb{F}F\left(X_{1},...,X_{k}\right),\left(Z_{1},...,Z_{k}\right)\right\rangle =\left.\frac{d}{ds}F\left(X_{1}+s\,Z_{1},...,X_{k}+s\,Z_{k}\right)\right\vert _{s=0}
\]
and
\[
\left\langle \mathbb{B}F\left(X_{1},...,X_{k}\right),Y\right\rangle =\left.\frac{d}{ds}F\left(u_{1}\left(s\right),...,u_{k}\left(s\right)\right)\right\vert _{s=0},
\]
respectively, where each $u_{i}:\left(-\varepsilon,\varepsilon\right)\rightarrow\mathcal{U}$
is a (horizontal)\ curve such that
\[
u_{i}\left(0\right)=X_{i},\ \ \ \left(\Pi\circ u_{i}\right)_{\ast}\left(\left.d/ds\right\vert _{0}\right)=Y\ \ \ \text{and}\ \ \ \frac{D}{Ds}u_{i}\left(s\right)=0.
\]

\end{remark}

As usual, by a tensor on $\mathcal{U}$ we mean a function $T:\mathcal{U}\times\cdot\cdot\cdot\times\mathcal{U}\rightarrow\mathbb{R}$,
on the Whitney sum of copies of $\mathcal{U}$, which is multi-linear
map when restricted to each fiber. When we write $T\left(X_{1},...,X_{k}\right)$,
it is implicit that all $X_{i}$'s are contained in the same fiber
of $\mathcal{U}$. 

Consider a tensor $\mathfrak{b}:T^{\ast}Q\times T^{\ast}Q\rightarrow\mathbb{R}$
and its related quadratic form $\mathfrak{q}:T^{\ast}Q\rightarrow\mathbb{R}:\alpha\longmapsto\mathfrak{b}\left(\alpha,\alpha\right).$
The following result is immediate.

\begin{proposition} For all $\alpha,\sigma\in T^{\ast}Q$,
\begin{equation}
\left\langle \mathbb{F}\mathfrak{b}\left(\alpha,\alpha\right),\left(\sigma,\sigma\right)\right\rangle =\left\langle \mathbb{F}\mathfrak{q}\left(\alpha\right),\sigma\right\rangle \ \ \ \text{and\ \ \ }\mathbb{B}\mathfrak{b}\left(\alpha,\alpha\right)=\mathbb{B}\mathfrak{q}\left(\alpha\right).\label{bq}
\end{equation}

\end{proposition}

\bigskip{}

Let $\omega$ be the canonical symplectic form on $T^{\ast}Q$.

\begin{proposition} Given $\alpha\in T^{\ast}Q$ and $V_{1},V_{2}\in T_{\alpha}T^{\ast}Q$,
and writing $\beta_{\alpha}\left(V_{1}\right)=v_{1}\oplus\sigma_{1}\ $and
$\beta_{\alpha}\left(V_{2}\right)=v_{2}\oplus\sigma_{2}$, we have
that 
\[
\omega\left(V_{1},V_{2}\right)=\left\langle \sigma_{2},v_{1}\right\rangle -\left\langle \sigma_{1},v_{2}\right\rangle +\left\langle \alpha,T\left(v_{1},v_{2}\right)\right\rangle ,
\]
being $T$ the torsion of $\nabla$. \end{proposition}

\emph{Proof. }Fixing a local chart $\left(U,\varphi\right)$ containing
$q=\pi\left(\alpha\right)$ and writing 
\[
\varphi^{\ast-1}\left(\alpha\right)=\left(\mathbf{q},\mathbf{p}\right)\;\;\text{and}\;\;\left(\varphi^{\ast-1}\right)_{\ast}\left(V_{\gamma}\right)=\left(\mathbf{q},\mathbf{p},\mathbf{\dot{q}_{\gamma}},\mathbf{\dot{p}_{\gamma}}\right),\;\;\;\gamma=1,2,
\]
it is well-known that $\omega\left(V_{1},V_{2}\right)=\dot{p}_{2,i}\,\dot{q}_{1}^{i}-\dot{p}_{1,i}\,\dot{q}_{2}^{i}$.
On the other hand, since
\[
\beta\left(\mathbf{q},\mathbf{p},\mathbf{\dot{q}_{\gamma}},\mathbf{\dot{p}_{\gamma}}\right)=\left(\mathbf{q},\mathbf{p}\right)\oplus\left(\mathbf{q},\mathbf{\dot{q}_{\gamma}}\right)\oplus\left(\mathbf{q},\mathbf{\dot{p}_{\gamma}}+\Gamma\left(\mathbf{q},\mathbf{p},\mathbf{\dot{q}_{\gamma}}\right)\right)
\]
{[}see $\left(\ref{bcpq}\right)${]}, we have that $\varphi_{\ast}\left(v_{\gamma}\right)=\left(\mathbf{q},\mathbf{\dot{q}_{\gamma}}\right)$
and $\varphi^{\ast-1}\left(\sigma_{\gamma}\right)=\left(\mathbf{q},\mathbf{\dot{p}_{\gamma}}+\Gamma\left(\mathbf{q},\mathbf{p},\mathbf{\dot{q}_{\gamma}}\right)\right)$,
and consequently
\[
\left\langle \sigma_{2},v_{1}\right\rangle -\left\langle \sigma_{1},v_{2}\right\rangle +\left\langle \alpha,T\left(v_{1},v_{2}\right)\right\rangle =\left(\dot{p}_{2,i}+\Gamma_{i}\left(\mathbf{q},\mathbf{p},\mathbf{\dot{q}}_{2}\right)\right)\,\dot{q}_{1}^{i}-\left(\dot{p}_{1,i}+\Gamma_{i}\left(\mathbf{q},\mathbf{p},\mathbf{\dot{q}}_{1}\right)\right)\,\dot{q}_{2}^{i}+p_{i}\,T^{i}\left(\mathbf{q},\mathbf{\dot{q}}_{1},\mathbf{\dot{q}}_{2}\right).
\]
So, we must show that $\Gamma_{i}\left(\mathbf{q},\mathbf{p},\mathbf{\dot{q}}_{2}\right)\,\dot{q}_{1}^{i}-\Gamma_{i}\left(\mathbf{q},\mathbf{p},\mathbf{\dot{q}}_{1}\right)\,\dot{q}_{2}^{i}=-p_{i}\,T^{i}\left(\mathbf{q},\mathbf{\dot{q}}_{1},\mathbf{\dot{q}}_{2}\right)$.
But 
\[
T^{i}\left(\mathbf{q},\mathbf{\dot{q}}_{1},\mathbf{\dot{q}}_{2}\right)=\Gamma_{kl}^{i}\left(\mathbf{q}\right)\ \dot{q}_{1}^{k}\,\dot{q}_{2}^{l}-\Gamma_{kl}^{i}\left(\mathbf{q}\right)\ \dot{q}_{2}^{k}\,\dot{q}_{1}^{l},
\]
from which the wanted result easily follows.\ \ \ $\square$

\bigskip{}

Assume that $\nabla$ is torsion-free, which we shall do from now
on. In terms of the diffeomorphisms $\beta$ and $\widetilde{\beta}$
we have the following result.

\begin{proposition} For all $v\in TQ$, $\alpha,\sigma\in T^{\ast}Q$,
on the same base point, 
\begin{equation}
\beta\circ\omega^{\sharp}\circ\widetilde{\beta}^{-1}\left(\alpha\oplus\sigma\oplus v\right)=\alpha\oplus v\oplus\left(-\sigma\right).\label{canb}
\end{equation}

\end{proposition}

\emph{Proof.} Following the notation of the previous proposition,
since $T=0$ (the torsion-free condition), we have that 
\[
\omega\left(\beta_{\alpha}^{-1}\left(v_{1}\oplus\sigma_{1}\right),\beta_{\alpha}^{-1}\left(v_{2}\oplus\sigma_{2}\right)\right)=\left\langle \sigma_{2},v_{1}\right\rangle -\left\langle \sigma_{1},v_{2}\right\rangle ,
\]
or equivalently
\[
\left\langle \left(\beta_{\alpha}^{\ast-1}\circ\omega^{\flat}\circ\beta_{\alpha}^{-1}\right)\left(v_{1}\oplus\sigma_{1}\right),\left(v_{2}\oplus\sigma_{2}\right)\right\rangle =\left\langle \sigma_{2},v_{1}\right\rangle -\left\langle \sigma_{1},v_{2}\right\rangle .
\]
This implies that $\left(\beta_{\alpha}^{\ast-1}\circ\omega^{\flat}\circ\beta_{\alpha}^{-1}\right)\left(v\oplus\sigma\right)=-\sigma\oplus v$
for all $v,\sigma$ on the base point of $\alpha$. Finally, using
the identity
\[
\left(\beta_{\alpha}^{\ast-1}\circ\omega^{\flat}\circ\beta_{\alpha}^{-1}\right)^{-1}=\beta_{\alpha}\circ\omega^{\sharp}\circ\beta_{\alpha}^{\ast},
\]
the proof is done.\ \ \ $\square$

\bigskip{}

Since the canonical Poisson bracket on $T^{\ast}Q$ is given by the
formula
\[
\{F,G\}(\alpha)=\left\langle dF(\alpha),\omega^{\sharp}\left(dG(\alpha)\right)\right\rangle ,\ \ \ \forall F,G\in C^{\infty}\left(T^{\ast}Q\right),
\]
using the last proposition and the Eq. $\left(\ref{bt}\right)$ we
easily arrive at the equation 
\begin{equation}
\{F,G\}(\alpha)=\left\langle \mathbb{B}F(\alpha),\mathbb{F}G(\alpha)\right\rangle -\left\langle \mathbb{B}G(\alpha),\mathbb{F}F(\alpha)\right\rangle .\label{fiberbasebracket}
\end{equation}
This identity will be central in the last section of the paper.

\bigskip{}

Finally, let us consider the next definition.

\begin{definition} \label{lyap}Let $P$ be a manifold and $X\in\mathfrak{X}(P)$
a vector field on $P$. Given a critical point $\alpha^{\bullet}\in P$
of $X$, i.e. $\alpha^{\bullet}$ such that $X\left(\alpha^{\bullet}\right)=0$,
a \textbf{Lyapunov function} for $X$ and $\alpha^{\bullet}$ is a
smooth function $\hat{H}:P\rightarrow\mathbb{R}$ satisfying
\begin{description}
\item [{L1}] $\hat{H}$ is positive definite w.r.t. $\alpha^{\bullet}$
(i.e. non-negative and null only at $\alpha^{\bullet}$);
\item [{L2}] $\left\langle d\hat{H}(\alpha),X(\alpha)\right\rangle \leq0$
for all $\alpha$. 
\end{description}
\end{definition}

As it is well-known, if such a function exists, then $\alpha^{\bullet}$
is a stable point. Moreover, if the inequality in \textbf{L2 }is strict
for all $\alpha\neq\alpha^{\bullet}$, then $\alpha^{\bullet}$ is
locally asymptotically stable, and if in addition $\hat{H}$ is a
proper function and $P$ is connected, then such a point is globally
asymptotically stable. For a proof of these results, see Ref. \cite{khalil}.

\section{Energy shaping method}

We present in this section the Hamiltonian \emph{side} of the energy
shaping method: the \emph{controlled Hamiltonians method} (as defined
in \cite{CBLMW02}), also known as the \emph{IDA-PBC method} \cite{idapbc}.
We describe a quite general version of the method, with its related
\emph{matching conditions}, in terms that are more convenient for
the present paper. For instance, we shall focus on Hamiltonian systems
on a cotangent bundle only. We shall progressively consider particular
situations to finally arrive at the case studied by Chang in Refs.
\cite{chang1}-\cite{chang4}, where particularly simple matching
conditions can be derived.

\subsection{The controlled Hamiltonians}

Fix a manifold $Q$, a function $H:T^{\ast}Q\rightarrow\mathbb{R}$
and a vertical subbundle $\mathcal{W}\subset\operatorname{ker}\pi_{\ast}\subset TT^{\ast}Q$
of the tangent bundle on $T^{\ast}Q$. Denote by $X_{H}:T^{\ast}Q\rightarrow TT^{\ast}Q$
the Hamiltonian vector field of $H$ w.r.t. the canonical symplectic
structure $\omega$ on $T^{\ast}Q$, i.e. $X_{H}:=\omega^{\sharp}\circ dH\in\mathfrak{X}\left(T^{\ast}Q\right)$.
Fix also a critical point $\alpha^{\bullet}\in T^{\ast}Q$ of $X_{H}$.
Note that the pair $\left(H,\mathcal{W}\right)$ defines an underactuated
Hamiltonian system on $Q$ (with Hamiltonian function $H$ and space
of actuators $\mathcal{W}$). It is clear that the rank of $\mathcal{W}$
represents the number of actuators. Suppose that we want to solve
the following problem.
\begin{description}
\item [{P.}] Find a control signal $Y\subset\mathcal{W}$, i.e. a vertical
vector field $Y\in\mathfrak{X}\left(T^{\ast}Q\right)$ with image
inside $\mathcal{W}$, such that the closed loop system defined by
$X_{H}+Y$ is stable at $\alpha^{\bullet}$. 
\end{description}
We shall call \emph{stabilization method }to any ``systematic procedure\textquotedblright \ that
enables us to solve the problem \textbf{P}. To be more precise, let
us consider the definitions below.

\begin{definition} \label{smd}Fix a manifold $Q$ and let $\mathfrak{U}$
be a subset of triples $\left(H,\mathcal{W},\alpha^{\bullet}\right)$,
where $\left(H,\mathcal{W}\right)$ is an underactuated system on
$Q$ and $\alpha^{\bullet}\in T^{\ast}Q$ is a critical point of $X_{H}$.
Given a triple $\left(H,\mathcal{W},\alpha^{\bullet}\right)\in\mathfrak{U}$,
denote by $\mathcal{S}_{H,\mathcal{W},\alpha^{\bullet}}\subset\mathfrak{X}\left(T^{\ast}Q\right)$
the subset of all the vector fields $Y\in\mathfrak{X}\left(T^{\ast}Q\right)$
solving \textbf{P}. We shall call \textbf{stabilization method} on
$\mathfrak{U}$ to any function\footnote{Note that a stabilization method on $\mathfrak{U}$ can also be seen
as relation in $\mathfrak{U}\times\mathfrak{X}\left(T^{\ast}Q\right)$.} $\digamma$ from $\mathfrak{U}$ to the power set of $\mathfrak{X}\left(T^{\ast}Q\right)$,
such that $\digamma\left(H,\mathcal{W},\alpha^{\bullet}\right)\subset\mathcal{S}_{H,\mathcal{W},\alpha^{\bullet}}$.
In addition, we shall say that $\digamma$ is \textbf{Lyapunov based}
if for each element $Y\in\digamma\left(H,\mathcal{W},\alpha^{\bullet}\right)$
a Lyapunov function for $X_{H}+Y$ and $\alpha^{\bullet}$ can be
exhibited (or at least exists).\footnote{The Massera's theorem \cite{mass} (and its various generalizations
-see \cite{kellett} for a review-) ensures that, if a smooth closed-loop
system is asymptotically stable, then a smooth Lyapunov function exists
for such a system. But the same can not be ensured if the system is
just stable (see \cite{Bacc}).} \end{definition}

\begin{definition} \label{ied}Given two stabilization methods $\digamma$
and $\digamma^{\prime}$ on the subsets $\mathfrak{U}$ and $\mathfrak{U}^{\prime}$,
respectively, we shall say that $\digamma$ is \textbf{included} in
$\digamma^{\prime}$ if
\[
\digamma\left(H,\mathcal{W},\alpha^{\bullet}\right)\subset\digamma^{\prime}\left(H,\mathcal{W},\alpha^{\bullet}\right),\ \ \ \forall\left(H,\mathcal{W},\alpha^{\bullet}\right)\in\mathfrak{U}\cap\mathfrak{U}^{\prime}.
\]
If both inclusions hold, we shall say that $\digamma$ and $\digamma^{\prime}$
are \textbf{equivalent} on $\mathfrak{U}\cap\mathfrak{U}^{\prime}$.
\end{definition}

\begin{remark} \label{generality} For methods on some common subset
$\mathfrak{U}$, the inclusion relation we just presented defines
a partial order. For such a partially ordered set, a maximal element
represents the most general way to systematically stabilize underactuated
systems in $\mathfrak{U}$.

The same can be said for the subset of Lyapunov based stabilization
methods. \end{remark}

The intention behind above definitions is to give a precise framework
to compare different methods and to establish what we mean by an equivalence
between them. Nonetheless, we will not construct the functions $\digamma$
when we describe the methods involved in this paper, but give a synthetic
explanation of the procedure they give rise to instead.

\bigskip{}

In the following, we shall define a Lyapunov based stabilization method
on the whole set of triples $\left(H,\mathcal{W},\alpha^{\bullet}\right)$,
known as the \emph{controlled Hamiltonians method}.

Assume that we are given a function $\hat{H}\in C^{\infty}\left(T^{\ast}Q\right)$,
an anti-symmetric tensor $B:T^{\ast}T^{\ast}Q\times T^{\ast}T^{\ast}Q\rightarrow\mathbb{R}$
(i.e. an almost-Poisson structure) on $T^{\ast}Q$ and two vertical
vector fields $Z_{g},Z_{d}\in\mathfrak{X}\left(T^{\ast}Q\right)$
such that:
\begin{enumerate}
\item $\hat{H}$ is positive-definite w.r.t. $\alpha^{\bullet}$,
\item $\left\langle d\hat{H}\left(\alpha\right),Z_{g}\left(\alpha\right)\right\rangle =0$
for all $\alpha\in T^{\ast}Q$, i.e. $Z_{g}$ is a \emph{gyroscopic}
force,
\item $B^{\sharp}\circ d\hat{H}+Z_{g}-X_{H}\subset\mathcal{W}$,
\item $\left\langle d\hat{H}\left(\alpha\right),Z_{d}\left(\alpha\right)\right\rangle \leq0$
for all $\alpha\in T^{\ast}Q$, i.e. $Z_{d}$ is a \emph{dissipative}
force,
\item $Z_{d}\subset\mathcal{W}$,
\item and $Z_{d}\left(\alpha^{\bullet}\right)=-Z_{g}\left(\alpha^{\bullet}\right)$. 
\end{enumerate}
Note that $1$ implies $d\hat{H}\left(\alpha^{\bullet}\right)=0$.
Then, because of $6$, defining $\hat{X}_{\hat{H}}:=B^{\sharp}\circ d\hat{H}$
and $\hat{X}:=\hat{X}_{\hat{H}}+Z_{g}+Z_{d}$, the point $\alpha^{\bullet}$
is a critical point of $\hat{X}$. Also, for all $\alpha\in T^{\ast}Q$,
\begin{equation}
\left\langle d\hat{H}\left(\alpha\right),\hat{X}\left(\alpha\right)\right\rangle =\left\langle d\hat{H}\left(\alpha\right),\hat{X}_{\hat{H}}\left(\alpha\right)+Z_{g}\left(\alpha\right)+Z_{d}\left(\alpha\right)\right\rangle =\left\langle d\hat{H}\left(\alpha\right),Z_{d}\left(\alpha\right)\right\rangle \leq0,\label{des}
\end{equation}
because of $2$, $4$ and the fact that $\left\langle d\hat{H}\left(\alpha\right),\hat{X}_{\hat{H}}\left(\alpha\right)\right\rangle =B\left(d\hat{H}\left(\alpha\right),d\hat{H}\left(\alpha\right)\right)=0$.
As a consequence, since $1$ coincides with \textbf{L1} and Eq. $\left(\ref{des}\right)$
coincides with \textbf{L2} (see Definition \ref{lyap}), $\hat{H}$
is a Lyapunov function for the dynamical system defined by $\hat{X}$
and the critical point $\alpha^{\bullet}$ . This says that such a
system is stable at $\alpha^{\bullet}$ (see Ref. \cite{khalil}).
So, defining 
\begin{equation}
Y:=\hat{X}-X_{H},\label{mes2}
\end{equation}
which belongs to $\mathcal{W}$ because of the points $3$ and $5$,
the problem \textbf{P }is solved. In particular, we have that 
\begin{equation}
\left\langle d\hat{H}(\alpha),X_{H}(\alpha)+Y(\alpha)\right\rangle \leq0.\label{des0}
\end{equation}
All that gives rise to the following procedure.

\begin{definition} \label{chm}Given an underactuated system $\left(H,\mathcal{W}\right)$
on $Q$ and a critical point $\alpha^{\bullet}\in T^{\ast}Q$ of $X_{H}$,
the \textbf{controlled Hamiltonians (CH) method} consists in finding
$\hat{H}$, $B$ and $Z_{g}$ satisfying $1$ and $2$ and solving
the equation {[}see point $3${]}
\begin{equation}
B^{\sharp}\circ d\hat{H}+Z_{g}-\omega^{\sharp}\circ dH\subset\mathcal{W};\label{mes11}
\end{equation}
finding $Z_{d}$ satisfying $4$, $5$ and $6$; and defining {[}see
$\left(\ref{mes2}\right)${]} 
\begin{equation}
Y:=B^{\sharp}\circ d\hat{H}+Z_{g}+Z_{d}-\omega^{\sharp}\circ dH.\label{mes3}
\end{equation}

\end{definition}

It is clear that above procedure defines, according to the Definition
\ref{smd}, a Lyapunov based stabilization method: its function $\digamma$
assigns to every triple $\left(H,\mathcal{W},\alpha^{\bullet}\right)$
a set of vector fields $\digamma\left(H,\mathcal{W},\alpha^{\bullet}\right)\subset\mathfrak{X}\left(T^{\ast}Q\right)$
given by Eq. $\left(\ref{mes3}\right)$ (and consequently solving
the problem \textbf{P}), where $\hat{H}$, $B$, $Z_{g}$ and $Z_{d}$
must fulfill the properties summarized in the last definition.

\begin{remark} The usual way of presenting the CH method is through
the idea of \emph{feedback equivalence} \cite{CBLMW02}. We shall
not explore this point of view here. \end{remark}

\subsection{A particular version}

The core of the CH method is Eq. $\left(\ref{mes11}\right)$, which
is a system of PDEs for $\hat{H}$, with unknown ``parameters\textquotedblright \ $B$
and $Z_{g}$. These PDEs are usually called \emph{matching conditions}.\footnote{The rest of the equations, the conditions $4$, $5$ and $6$, define
algebraic conditions for $Z_{d}$ that will be studied later.} Different assumptions on the original underactuated system $\left(H,\mathcal{W}\right)$,
and particular ansatzs for the unknowns $\hat{H}$, $B$ and $Z_{g}$,
give rise to particular forms of $\left(\ref{mes11}\right)$ and,
consequently, to particular versions of the method. (In terms of Definition
\ref{ied}, we have in this way different included methods.) For instance,
let us assume that $H,\hat{H}:T^{\ast}Q\rightarrow\mathbb{R}$ are
hyper-regular, i.e. $\mathbb{F}H,\mathbb{F}\hat{H}:T^{\ast}Q\rightarrow TQ$
are linear bundle isomorphisms {[}see $\left(\ref{fiberder}\right)${]},
and also
\begin{equation}
\mathbb{F}H=\mathbb{F}H^{\ast}\ \ \ \text{and\ \ \ }\mathbb{F}\hat{H}=\mathbb{F}\hat{H}^{\ast},\label{symm}
\end{equation}
i.e. their fiber derivatives are symmetric. In addition, fix a torsion-free
connection on $T^{\ast}Q$ and assume that {[}recall Eqs. \eqref{beta}
and \eqref{betam}{]} 
\begin{equation}
\beta\circ B^{\sharp}\circ\widetilde{\beta}^{-1}\left(\alpha\oplus\sigma\oplus v\right)=\alpha\oplus\Psi\left(v\right)\oplus\Psi^{\ast}\left(-\sigma\right),\label{wc}
\end{equation}
for some fiber bundle morphism $\Psi:TQ\rightarrow TQ$ {[}compare
to Eq. $\left(\ref{canb}\right)${]}. 

\begin{remark}Note that, according to \eqref{Yconbeta}, each $\mathcal{W}_{\alpha}\subset T_{\alpha}T^{\ast}Q$
is defined by the unique subspace $W_{\alpha}\subset T_{\pi\left(\alpha\right)}^{\ast}Q$
such that
\begin{equation}
\mathcal{W}_{\alpha}=\operatorname{vlift}_{\alpha}\left(W_{\alpha}\right)=\beta^{-1}\left(\alpha\oplus0\oplus W_{\alpha}\right).\label{wa}
\end{equation}
Also {[}see Eq. \eqref{fieldconbeta}{]}, the vertical vector field
$Z_{g}$ can be written 
\begin{equation}
Z_{g}\left(\alpha\right)=\beta^{-1}\left(\alpha\oplus0\oplus z_{g}\left(\alpha\right)\right),\label{zg}
\end{equation}
for a unique fiber preserving map $z_{g}:T^{\ast}Q\rightarrow T^{\ast}Q$.
(Idem $Z_{d}$.) \end{remark}

\begin{proposition}Under above assumptions and notation, the matching
conditions $\left(\ref{mes11}\right)$ reduce to {[}see Eq. \eqref{zg}{]}
\begin{equation}
\left\langle \mathbb{B}\hat{H}\left(\alpha\right),\mathbb{F}H\left(\sigma\right)\right\rangle -\left\langle \mathbb{B}H\left(\alpha\right),\mathbb{F}\hat{H}\left(\sigma\right)\right\rangle -\left\langle z_{g}\left(\alpha\right),\mathbb{F}\hat{H}\left(\sigma\right)\right\rangle =0,\:\:\:\forall\sigma\in\hat{W}_{\alpha},\label{ida1}
\end{equation}
where {[}see Eq. \eqref{wa}{]}
\begin{equation}
\hat{W}_{\alpha}:=\left(\mathbb{F}\hat{H}\left(W_{\alpha}\right)\right)^{0}=\mathbb{F}\hat{H}^{-1}\left(W_{\alpha}^{0}\right),\:\:\:\forall\alpha\in T^{\ast}Q.\label{wac}
\end{equation}
In particular, the unknowns are $\hat{H}$ and $z_{g}$ \textbf{only}.
\end{proposition}

\emph{Proof.} From Eqs. $\left(\ref{bt}\right)$ and \eqref{wc},
we have that
\begin{align}
\beta\left(\hat{X}_{\hat{H}}\left(\alpha\right)\right) & =\beta\circ B^{\sharp}\left(d\hat{H}\left(\alpha\right)\right)=\beta\circ B^{\sharp}\circ\widetilde{\beta}^{-1}\left(\widetilde{\beta}\left(d\hat{H}\left(\alpha\right)\right)\right)=\beta\circ B^{\sharp}\circ\widetilde{\beta}^{-1}\left(\alpha\oplus\mathbb{B}\hat{H}\left(\alpha\right)\oplus\mathbb{F}\hat{H}\left(\alpha\right)\right)\nonumber \\
 & =\alpha\oplus\Psi\left(\mathbb{F}\hat{H}\left(\alpha\right)\right)\oplus\Psi^{\ast}\left(-\mathbb{B}\hat{H}\left(\alpha\right)\right).\label{bxc}
\end{align}
Similarly, from Eqs. $\left(\ref{bt}\right)$ and $\left(\ref{canb}\right)$,
\begin{equation}
\beta\left(X_{H}\left(\alpha\right)\right)=\alpha\oplus\mathbb{F}H\left(\alpha\right)\oplus\left(-\mathbb{B}H\left(\alpha\right)\right).\label{bx}
\end{equation}
As a consequence, using Equations $\left(\ref{wa}\right)$, $\left(\ref{zg}\right)$,
$\left(\ref{bxc}\right)$ and $\left(\ref{bx}\right)$, it easily
follows that $\left(\ref{mes11}\right)$ reduces to
\[
\Psi\left(\mathbb{F}\hat{H}\left(\alpha\right)\right)=\mathbb{F}H\left(\alpha\right)\ \ \ \ \text{and\ \ \ }-\Psi^{\ast}\left(\mathbb{B}\hat{H}\left(\alpha\right)\right)+z_{g}\left(\alpha\right)+\mathbb{B}H\left(\alpha\right)\in W_{\alpha},
\]
for all $\alpha\in T^{\ast}Q$. This implies that $\Psi=\mathbb{F}H\circ\mathbb{F}\hat{H}^{-1}$
and, taking Eq. $\left(\ref{symm}\right)$ into account,
\[
-\mathbb{F}\hat{H}^{-1}\circ\mathbb{F}H\circ\mathbb{B}\hat{H}\left(\alpha\right)+z_{g}\left(\alpha\right)+\mathbb{B}H\left(\alpha\right)\in W_{\alpha}.
\]
It only rests to use the Eq. \eqref{wac} in order to end the proof.\ \ \ $\square$

\bigskip{}

\begin{remark} \label{cpc1}It is clear from $\left(\ref{bx}\right)$
that $\alpha^{\bullet}$ is a critical point of $X_{H}$ if and only
if $\mathbb{F}H\left(\alpha^{\bullet}\right)=0$ and $\mathbb{B}H\left(\alpha^{\bullet}\right)=0$.
\end{remark}

Let us mention that, according to Eqs. $\left(\ref{fieldconbeta}\right)$
and $\left(\ref{mes3}\right)$, each control law $Y$ of the method
is now given by the vertical lift of the fiber preserving map
\begin{equation}
y\left(\alpha\right):=z_{d}\left(\alpha\right)+z_{g}\left(\alpha\right)-\mathbb{F}\hat{H}^{-1}\circ\left(\mathbb{F}H\circ\mathbb{B}\hat{H}\left(\alpha\right)-\mathbb{F}\hat{H}\circ\mathbb{B}H\left(\alpha\right)\right).\label{y0}
\end{equation}
Also, 
\begin{equation}
\left\langle z_{g}\left(\alpha\right),\mathbb{F}\hat{H}\left(\alpha\right)\right\rangle =0\label{zg0}
\end{equation}
from point $2$ above (the gyroscopic condition),
\begin{equation}
\left\langle z_{d}\left(\alpha\right),\mathbb{F}\hat{H}\left(\alpha\right)\right\rangle \leq0\ \ \ \ \text{and\ \ \ \ }z_{d}\left(\alpha\right)\in W_{\alpha}\label{Zd}
\end{equation}
according to points $4$ and $5$, and 
\begin{equation}
z_{d}\left(\alpha^{\bullet}\right)=-z_{g}\left(\alpha^{\bullet}\right)\label{zdp}
\end{equation}
from point $6$.

\subsection{The \emph{kinetic} and \emph{potential} matching conditions}

\label{kpmc}

Let us further restrict the original underactuated system $\left(H,\mathcal{W}\right)$
and the unknowns $\hat{H}$ and $z_{g}$ of $\left(\ref{ida1}\right)$.
Assume first that $H:T^{\ast}Q\rightarrow\mathbb{R}$ is a \emph{simple}
Hamiltonian function, i.e. 
\[
H=\mathfrak{H}+h\circ\pi\ \ \ \text{with\ \ \ }\mathfrak{H}(\alpha):=\frac{1}{2}\left\langle \alpha,\rho^{\sharp}(\alpha)\right\rangle ,\ \ \forall\alpha\in T^{\ast}Q,
\]
where $h\in C^{\infty}(Q)$ and $\rho$ is a Riemannian metric on
$Q$. The first and second terms in $H$ are known as the \emph{kinetic}
and \emph{potential} terms. Note that $\mathfrak{H}$ is the quadratic
form of the tensor 
\[
\mathfrak{b}:T^{\ast}Q\times T^{\ast}Q\rightarrow\mathbb{R}:\left(\alpha,\sigma\right)\longmapsto\frac{1}{2}\left\langle \alpha,\rho^{\sharp}(\sigma)\right\rangle .
\]
Also note that $\mathbb{F}H=\mathbb{F}\mathfrak{H}=\rho^{\sharp}$
{[}for the first equality, recall the Eq. $\left(\ref{basic}\right)$
of Proposition \ref{bas}{]}. This implies that $\mathbb{F}H$ is
a symmetric linear bundle isomorphism. Choosing a coordinate chart
$\left(U,\varphi\right)$ on $Q$ and its induced one on $T^{\ast}Q$
{[}see \eqref{icor} and \eqref{icor1}{]}, i.e. choosing canonical
coordinates, and denoting by $\mathbb{H}\left(\mathbf{q}\right)$
the coordinate matrix representation at $q$ of the Riemannian metric
$\rho$, we can write 
\begin{equation}
H(\mathbf{q},\mathbf{p})=\frac{1}{2}p_{i}\,\mathbb{H}^{ij}(\mathbf{q})\,p_{j}+h(\mathbf{q}).\label{lhe}
\end{equation}
Of course, the symmetric condition $\mathbb{F}H=\mathbb{F}H^{\ast}$
translates to $\mathbb{H}$ as $\mathbb{H}^{ij}=\mathbb{H}^{ji}$.

\begin{remark} \label{cps}As it is well-known, for Hamiltonian systems
defined by a simple function, the critical points are of the form
$\alpha^{\bullet}=\left(q^{\bullet},0\right)$, with $dh\left(q^{\bullet}\right)=0$
(use Proposition \ref{bas} and Remark \ref{cpc1}). \end{remark}

Regarding the unknowns of $\left(\ref{ida1}\right)$, assume that
$\hat{H}$ is simple too. We shall use for $\hat{H}$ an analogue
notation to that we used for $H$. For instance, we shall write $\hat{H}=\hat{\mathfrak{H}}+\hat{h}\circ\pi$.
Thus, $\mathbb{F}\hat{H}=\mathbb{F}\hat{\mathfrak{H}}$. Also, assume
that $z_{g}:T^{\ast}Q\rightarrow T^{\ast}Q$ is given by
\begin{equation}
\left\langle z_{g}\left(\alpha\right),v\right\rangle =\mathfrak{Z}_{g}\left(\alpha,\alpha,\mathbb{F}\hat{\mathfrak{H}}^{-1}\left(v\right)\right),\ \ \ \forall q\in Q,\ \alpha\in T_{q}^{\ast}Q,\ v\in T_{q}Q,\label{zZg}
\end{equation}
for some tensor field $\mathfrak{Z}_{g}:T^{\ast}Q\times T^{\ast}Q\times T^{\ast}Q\rightarrow\mathbb{R}$.
This particular choice for the map $z_{g}$ implies that it is quadratic
in $\alpha$. Note that the tensor field $\mathfrak{Z}_{g}$ can be
assumed symmetric in its first two arguments, i.e. 
\begin{equation}
\mathfrak{Z}_{g}\left(\alpha_{1},\alpha_{2},\alpha\right)=\mathfrak{Z}_{g}\left(\alpha_{2},\alpha_{1},\alpha\right).\label{g1}
\end{equation}
In addition, using the gyroscopic condition {[}see Eq. $\left(\ref{zg0}\right)${]},
it is clear that 
\begin{equation}
\mathfrak{Z}_{g}\left(\alpha,\alpha,\alpha\right)=0,\label{g2}
\end{equation}
for all $\alpha\in T_{q}^{\ast}Q$. Reciprocally, any tensor satisfying
\eqref{g2} gives rise, through $\left(\ref{zg}\right)$ and $\left(\ref{zZg}\right)$,
to a gyroscopic force.

\bigskip{}

Coming back to the matching conditions, since the kinetic terms of
$H$ and $\hat{H}$ are quadratic functions and their potential terms
are basic functions, the next result can be easily proved.

\begin{proposition}Under above assumptions and notation, the matching
conditions $\left(\ref{ida1}\right)$ decompose into two equations:
\begin{equation}
\left\langle \mathbb{B}\hat{\mathfrak{H}}\left(\alpha\right),\mathbb{F}\mathfrak{H}\left(\sigma\right)\right\rangle -\left\langle \mathbb{B}\mathfrak{H}\left(\alpha\right),\mathbb{F}\hat{\mathfrak{H}}\left(\sigma\right)\right\rangle -\mathfrak{Z}_{g}\left(\alpha,\alpha,\sigma\right)=0,\label{kei}
\end{equation}
the \emph{kinetic matching conditions}, and 
\begin{equation}
\left\langle d\hat{h}\left(\pi\left(\sigma\right)\right),\mathbb{F}\mathfrak{H}\left(\sigma\right)\right\rangle -\left\langle dh\left(\pi\left(\sigma\right)\right),\mathbb{F}\hat{\mathfrak{H}}\left(\sigma\right)\right\rangle =0,\label{pei}
\end{equation}
the \emph{potential matching conditions}. They must be satisfied for
all 
\begin{equation}
\alpha\in T^{\ast}Q\ \ \ \text{and\ \ \ }\sigma\in\hat{W}_{\alpha}.\label{sa}
\end{equation}

\end{proposition}

\begin{remark} \label{ort}Note that {[}see Eq. $\left(\ref{wac}\right)${]}
$\hat{W}_{\alpha}=\mathbb{F}\hat{\mathfrak{H}}^{-1}\left(W_{\alpha}^{0}\right)$
is the orthogonal complement of $W_{\alpha}$ w.r.t. the bilinear
$\hat{\mathfrak{b}}$.

\end{remark}

In local coordinates $\left(U,\varphi\right)$, combining $\left(\ref{lfd}\right)$,
$\left(\ref{lbd}\right)$ and $\left(\ref{lhe}\right)$, 
\begin{equation}
\begin{aligned}\left\langle \mathbb{B}\hat{\mathfrak{H}}\left(\alpha\right),\mathbb{F}\mathfrak{H}\left(\sigma\right)\right\rangle -\left\langle \mathbb{B}\mathfrak{H}\left(\alpha\right),\mathbb{F}\hat{\mathfrak{H}}\left(\sigma\right)\right\rangle  & =\frac{1}{2}\left(\frac{\partial\hat{\mathbb{H}}^{ij}(\mathbf{q})}{\partial q^{k}}\,\mathbb{H}^{kl}(\mathbf{q})-\frac{\partial\mathbb{H}^{ij}(\mathbf{q})}{\partial q^{k}}\,\hat{\mathbb{H}}^{kl}(\mathbf{q})\right)\,p_{i}p_{j}\widetilde{p}_{l}\\
 & +\Gamma_{ks}^{j}(\mathbf{q})\,\left(\hat{\mathbb{H}}^{is}(\mathbf{q})\,\mathbb{H}^{kl}(\mathbf{q})-\mathbb{H}^{is}(\mathbf{q})\,\hat{\mathbb{H}}^{kl}(\mathbf{q})\right)\,p_{i}p_{j}\widetilde{p}_{l}
\end{aligned}
\label{id1}
\end{equation}
and
\begin{equation}
\left\langle d\hat{h}\left(\pi\left(\sigma\right)\right),\mathbb{F}\mathfrak{H}\left(\sigma\right)\right\rangle -\left\langle dh\left(\pi\left(\sigma\right)\right),\mathbb{F}\hat{\mathfrak{H}}\left(\sigma\right)\right\rangle =\left(\frac{\partial\hat{h}(\mathbf{q})}{\partial q^{k}}\,\mathbb{H}^{kl}(\mathbf{q})-\frac{\partial h(\mathbf{q})}{\partial q^{k}}\,\hat{\mathbb{H}}^{kl}(\mathbf{q})\right)\,\widetilde{p}_{l},\label{id2}
\end{equation}
for $\alpha=\varphi^{\ast}\left(\mathbf{q},\mathbf{p}\right)$ and
$\sigma=\varphi^{\ast}\left(\mathbf{q},\mathbf{\widetilde{p}}\right)$.
Thus, Eqs. $\left(\ref{kei}\right)$ and $\left(\ref{pei}\right)$
translate to
\begin{equation}
\begin{aligned}\left(\frac{\partial\hat{\mathbb{H}}^{ij}(\mathbf{q})}{\partial q^{k}}\,\mathbb{H}^{kl}(\mathbf{q})-\frac{\partial\mathbb{H}^{ij}(\mathbf{q})}{\partial q^{k}}\,\hat{\mathbb{H}}^{kl}(\mathbf{q})-2\,\mathfrak{Z}_{g}^{ijl}(\mathbf{q})\right)\,p_{i}p_{j}\widetilde{p}_{l}\\
+2\,\Gamma_{ks}^{j}(\mathbf{q})\,\left(\hat{\mathbb{H}}^{is}(\mathbf{q})\,\mathbb{H}^{kl}(\mathbf{q})-\mathbb{H}^{is}(\mathbf{q})\,\hat{\mathbb{H}}^{kl}(\mathbf{q})\right)\,p_{i}p_{j}\widetilde{p}_{l} & =0
\end{aligned}
\label{ke0}
\end{equation}
and 
\[
\left(\frac{\partial\hat{h}(\mathbf{q})}{\partial q^{k}}\,\mathbb{H}^{kl}(\mathbf{q})-\frac{\partial h(\mathbf{q})}{\partial q^{k}}\,\hat{\mathbb{H}}^{kl}(\mathbf{q})\right)\,\widetilde{p}_{l}=0,
\]
for all $q\in U$, $\mathbf{p}\in\left(\varphi_{q}^{\ast}\right)^{-1}\left(T_{q}^{\ast}Q\right)$
and $\mathbf{\widetilde{p}}\in\left(\varphi_{q}^{\ast}\right)^{-1}\left(\hat{W}_{\varphi_{q}^{\ast}\left(\mathbf{p}\right)}\right)$.
(Of course, the numbers $\mathfrak{Z}_{g}^{ijl}(\mathbf{q})$ are
the coefficients of the coordinate matrix representation of $\mathfrak{Z}_{g}$
at $q$). 

\bigskip{}

Equations $\left(\ref{kei}\right)$, $\left(\ref{pei}\right)$ and
$\left(\ref{sa}\right)$ define a generalized version of the traditional
IDA-PBC method \cite{idapbc}, also known as the \emph{energy plus
force shaping method}. In the following subsection, we shall make
two more assumptions that will drive us to another set of matching
conditions, originally\footnote{Actually, such matching conditions first appeared in \cite{ham},
but without a derivation.} studied by Chang in \cite{chang3,chang4} (see also \cite{chang1,chang2}
for the Lagrangian counterpart).

\bigskip{}

To end this section, note that under all above assumptions {[}see
Eq. $\left(\ref{y0}\right)${]}
\begin{equation}
\begin{aligned}y\left(\alpha\right) & =z_{d}\left(\alpha\right)+\mathfrak{Z}_{g}\left(\alpha,\alpha,\mathbb{F}\hat{H}^{-1}\left(\cdot\right)\right)-\mathbb{F}\hat{H}^{-1}\circ\left(\mathbb{F}H\circ\mathbb{B}\hat{H}\left(\alpha\right)-\mathbb{F}\hat{H}\circ\mathbb{B}H\left(\alpha\right)\right)\\
 & =z_{d}\left(\alpha\right)+\mathfrak{Z}_{g}\left(\alpha,\alpha,\mathbb{F}\hat{\mathfrak{H}}^{-1}\left(\cdot\right)\right)-\mathbb{F}\hat{\mathfrak{H}}^{-1}\circ\left(\mathbb{F}\mathfrak{H}\circ\mathbb{B}\hat{\mathfrak{H}}\left(\alpha\right)-\mathbb{F}\hat{\mathfrak{H}}\circ\mathbb{B}\mathfrak{H}\left(\alpha\right)\right)\\
 & -\mathbb{F}\hat{\mathfrak{H}}^{-1}\circ\left(\mathbb{F}\mathfrak{H}\circ d\hat{h}\left(\pi\left(\alpha\right)\right)-\mathbb{F}\hat{\mathfrak{H}}\circ dh\left(\pi\left(\alpha\right)\right)\right).
\end{aligned}
\label{y}
\end{equation}

\subsection{\emph{Simple} matching conditions}

\label{cmc}

First, assume that there exists a subbundle $W\subset T^{\ast}Q$
such that {[}recall Eq. $\left(\ref{wa}\right)${]} 
\begin{equation}
\mathcal{W}_{\alpha}=\operatorname{vlift}_{\alpha}\left(W_{\pi\left(\alpha\right)}\right),\ \ \ \forall\alpha\in T^{\ast}Q.\label{ww}
\end{equation}
Following $\left(\ref{wac}\right)$, let us define another subbundle
of $T^{\ast}Q$,
\begin{equation}
\hat{W}:=\left(\mathbb{F}\hat{\mathfrak{H}}\left(W\right)\right)^{0}=\mathbb{F}\hat{\mathfrak{H}}^{-1}\left(W^{0}\right).\label{wac2}
\end{equation}

\begin{remark} \label{ort1}Note that, according to Remark \ref{ort},
$\hat{W}$ is the orthogonal complement of $W$ w.r.t. to the tensor
$\hat{\mathfrak{b}}$. In particular, we have that $T^{\ast}Q=W\oplus\hat{W}$.
\end{remark}

Under these assumptions, Chang showed in \cite{chang1} that there
exists a solution $\left(\hat{\mathfrak{H}},\mathfrak{Z}_{g}\right)$
of $\left(\ref{kei}\right)$ if and only if $\hat{\mathfrak{H}}$
satisfies
\[
\left\langle \mathbb{B}\hat{\mathfrak{H}}\left(\sigma\right),\mathbb{F}\mathfrak{H}\left(\sigma\right)\right\rangle -\left\langle \mathbb{B}\mathfrak{H}\left(\sigma\right),\mathbb{F}\hat{\mathfrak{H}}\left(\sigma\right)\right\rangle =0,\ \ \ \forall\sigma\in\hat{W}.
\]
Moreover, it can be shown by using elementary tensor algebra that:

\begin{proposition}$\left(\hat{\mathfrak{H}},\mathfrak{Z}_{g}\right)$
is a solution of $\left(\ref{kei}\right)$ if and only if $\hat{\mathfrak{H}}$
satisfies the above equation and $\mathfrak{Z}_{g}$ is given as follows:
\begin{enumerate}
\item define $\Upsilon:T^{\ast}Q\times T^{\ast}Q\times T^{\ast}Q\rightarrow\mathbb{R}$
as 
\begin{equation}
\Upsilon\left(\alpha_{1},\alpha_{2},\alpha_{3}\right):=\left\langle \mathbb{B}\hat{\mathfrak{b}}\left(\alpha_{1},\alpha_{2}\right),\mathbb{F}\mathfrak{H}\left(\alpha_{3}\right)\right\rangle -\left\langle \mathbb{B}\mathfrak{b}\left(\alpha_{1},\alpha_{2}\right),\mathbb{F}\hat{\mathfrak{H}}\left(\alpha_{3}\right)\right\rangle ;\label{T}
\end{equation}
\item fix a tensor $A:W\times W\times\hat{W}\rightarrow\mathbb{R}$ satisfying
\[
A\left(\alpha_{1},\alpha_{2},\alpha\right)=-A\left(\alpha_{2},\alpha_{1},\alpha\right),\ \ \ \forall\alpha_{1},\alpha_{2}\in W,\ \alpha\in\hat{W},
\]
and a tensor $B:W\times W\times W\rightarrow\mathbb{R}$ satisfying
$\left(\ref{g1}\right)$ and $\left(\ref{g2}\right)$ along $W$;
\item and finally define $\mathfrak{Z}_{g}:T^{\ast}Q\times T^{\ast}Q\times T^{\ast}Q\rightarrow\mathbb{R}$
as\footnote{Chang made a similar construction to show the existence of solutions
of $\left(\ref{kei}\right)$ (see Eqs. $2.8$ to $2.11$ in Ref. \cite{chang4},
and replace $\Upsilon$ and $\mathfrak{Z}_{g}$ by $S$ and $C$,
respectively), but taking $A=0$ and $B=0$.}
\begin{equation}
\left\{ \begin{array}{l}
\mathfrak{Z}_{g}\left(\alpha_{1},\alpha_{2},\sigma\right):=\Upsilon\left(\alpha_{1},\alpha_{2},\sigma\right),\ \ \ \mathfrak{Z}_{g}\left(\sigma_{1},\sigma_{2},\gamma\right):=-\Upsilon\left(\gamma,\sigma_{2},\sigma_{1}\right)-\Upsilon\left(\gamma,\sigma_{1},\sigma_{2}\right),\\
\mathfrak{Z}_{g}\left(\gamma_{1},\sigma,\gamma_{2}\right):=\mathfrak{Z}_{g}\left(\sigma,\gamma_{1},\gamma_{2}\right):=-\frac{1}{2}\left[\Upsilon\left(\gamma_{1},\gamma_{2},\sigma\right)+A\left(\gamma_{1},\gamma_{2},\sigma\right)\right],\\
\mathfrak{Z}_{g}\left(\gamma_{1},\gamma_{2},\gamma_{3}\right):=B\left(\gamma_{1},\gamma_{2},\gamma_{3}\right),
\end{array}\right.\label{Zg}
\end{equation}
with $\alpha_{i}\in T^{\ast}Q$, $\sigma,\sigma_{i}\in\hat{W}$ and
$\gamma,\gamma_{i}\in W$. 
\end{enumerate}
\end{proposition}

Then, a solution $\left(\hat{\mathfrak{H}},\hat{h},\mathfrak{Z}_{g}\right)$
of $\left(\ref{kei}\right)$, $\left(\ref{pei}\right)$ and $\left(\ref{sa}\right)$
can be found if and only if we solve the equations
\begin{align}
\left\langle \mathbb{B}\hat{\mathfrak{H}}\left(\sigma\right),\mathbb{F}\mathfrak{H}\left(\sigma\right)\right\rangle -\left\langle \mathbb{B}\mathfrak{H}\left(\sigma\right),\mathbb{F}\hat{\mathfrak{H}}\left(\sigma\right)\right\rangle  & =0,\label{kkechi}\\
\left\langle d\hat{h}\left(\pi\left(\sigma\right)\right),\mathbb{F}\mathfrak{H}\left(\sigma\right)\right\rangle -\left\langle dh\left(\pi\left(\sigma\right)\right),\mathbb{F}\hat{\mathfrak{H}}\left(\sigma\right)\right\rangle  & =0,\ \ \ \forall\sigma\in\hat{W},\label{pechi}
\end{align}
for $\hat{\mathfrak{H}}$ and $\hat{h}$ \textbf{only}. These are
the new matching conditions that we mentioned above, which we shall
call the \emph{Chang's }or \emph{simple matching conditions}. The
local counterpart reads
\begin{equation}
\left(\frac{\partial\hat{\mathbb{H}}^{ij}(\mathbf{q})}{\partial q^{k}}\,\mathbb{H}^{kl}(\mathbf{q})-\frac{\partial\mathbb{H}^{ij}(\mathbf{q})}{\partial q^{k}}\,\hat{\mathbb{H}}^{kl}(\mathbf{q})\right)\,p_{i}p_{j}p_{l}=0\label{kec}
\end{equation}
{[}see Eq. $\left(\ref{id1}\right)$ for $\alpha=\sigma$ and use
the torsion-free condition{]} and 
\begin{equation}
\left(\frac{\partial\hat{h}(\mathbf{q})}{\partial q^{k}}\,\mathbb{H}^{kl}(\mathbf{q})-\frac{\partial h(\mathbf{q})}{\partial q^{k}}\,\hat{\mathbb{H}}^{kl}(\mathbf{q})\right)\,p_{l}=0\label{pec}
\end{equation}
{[}see Eq. $\left(\ref{id2}\right)${]}, for all $q\in U$ and $\mathbf{p}\in\left(\varphi_{q}^{\ast}\right)^{-1}\left(\hat{W}_{q}\right)$,
or equivalently {[}see $\left(\ref{wac2}\right)${]} $\varphi_{\ast,q}^{-1}\left(\hat{\mathbb{H}}(\mathbf{q})\cdot\mathbf{p}\right)\in W_{q}^{0}$.

\begin{remark} Above equations are (up to a sign) Eqs. $\left(2.21\right)$
and $\left(2.29\right)$ of \cite{chang4} for 
\[
M^{ij}=\mathbb{H}^{ij},\ \ \ \widehat{M}^{ij}=\hat{\mathbb{H}}^{ij},\ \ \ V=h,\ \ \ \widehat{V}=\hat{h},\ \ \text{and\ \ \ }G^{\perp}=W^{0}.
\]

\end{remark}

Now, let us study the Eqs. \eqref{Zd} and \eqref{zdp} in the present
situation.

\begin{remark} If $\alpha^{\bullet}=\left(q^{\bullet},0\right)$
is a critical point for $X_{H}$, since $\mathfrak{Z}_{g}\left(0,0,\mathbb{F}\hat{\mathfrak{H}}^{-1}\left(\cdot\right)\right)=0$,
condition $\left(\ref{zdp}\right)$ for $z_{d}$ reduces to $z_{d}\left(\alpha^{\bullet}\right)=0$.
\end{remark}

According to the last remark and using that
\[
\left\langle z_{d}\left(\alpha\right),\mathbb{F}\hat{H}\left(\alpha\right)\right\rangle =\left\langle z_{d}\left(\alpha\right),\mathbb{F}\hat{\mathfrak{H}}\left(\alpha\right)\right\rangle =\hat{\mathfrak{b}}\left(z_{d}\left(\alpha\right),\alpha\right),
\]
such equations can be written
\begin{equation}
\hat{\mathfrak{b}}\left(z_{d}\left(\alpha\right),\alpha\right)\leq0,\text{\ \ \ }z_{d}\left(\alpha\right)\in W\,\ \ \ \text{and}\ \ \ z_{d}\left(\alpha^{\bullet}\right)=0.\label{szd}
\end{equation}
Suppose that we have a solution $z_{d}$ of $\left(\ref{szd}\right)$,
and consider the orthogonal projection $P$ with image $W$ (see Remark
\ref{ort1}). Then, 
\[
\hat{\mathfrak{b}}\left(z_{d}\left(\alpha\right),\alpha\right)=\hat{\mathfrak{b}}\left(z_{d}\left(\alpha\right),P\left(\alpha\right)\right)=-\mu\left(\alpha\right)
\]
for some non-negative function $\mu:T^{\ast}Q\rightarrow\mathbb{R}$
such that $\mu\left(\sigma\right)=0$ for all $\sigma\in\hat{W}$.
Fixing $\alpha_{0}\notin\hat{W}$, i.e. $P\left(\alpha_{0}\right)\neq0$,
and using elementary linear algebra, it follows that
\[
z_{d}\left(\alpha_{0}\right)=x_{0}-\frac{\hat{\mathfrak{b}}\left(x_{0},P\left(\alpha_{0}\right)\right)+\mu\left(\alpha_{0}\right)}{\hat{\mathfrak{b}}\left(P\left(\alpha_{0}\right),P\left(\alpha_{0}\right)\right)}\ P\left(\alpha_{0}\right),\ \ \ \text{for some \ \ \ }x_{0}\in W.
\]
In addition, since the complementary subset of $\hat{W}$ in $T^{\ast}Q$
is an open dense submanifold, this means that $z_{d}$ must be given
by the expression
\begin{equation}
x\left(\alpha\right)-\frac{\hat{\mathfrak{b}}\left(x\left(\alpha\right),P\left(\alpha\right)\right)+\mu\left(\alpha\right)}{\hat{\mathfrak{b}}\left(P\left(\alpha\right),P\left(\alpha\right)\right)}\ P\left(\alpha\right),\ \ \ \forall\alpha\notin\hat{W},\label{cif}
\end{equation}
for some fiber preserving map $x:T^{\ast}Q\rightarrow T^{\ast}Q$
with image inside $W$. And since $z_{d}\left(\alpha^{\bullet}\right)=0$,
\begin{equation}
\underset{\alpha\rightarrow\alpha^{\bullet}}{\operatorname{Lim}}\left(x\left(\alpha\right)-\frac{\hat{\mathfrak{b}}\left(x\left(\alpha\right),P\left(\alpha\right)\right)+\mu\left(\alpha\right)}{\hat{\mathfrak{b}}\left(P\left(\alpha\right),P\left(\alpha\right)\right)}\ P\left(\alpha\right)\right)=0.\label{lim}
\end{equation}
Reciprocally, it is easy to see that, if $z_{d}$ is a smooth map
given by the formula $\left(\ref{cif}\right)$ and satisfies $\left(\ref{lim}\right)$,
then it satisfies $\left(\ref{szd}\right)$. Concluding,

\begin{proposition}\label{stepsiii}If we:
\begin{description}
\item [{i.}] fix a non-negative function $\mu:T^{\ast}Q\rightarrow\mathbb{R}$
such that $\mu\left(\sigma\right)=0$ for all $\sigma\in\hat{W}$;
\item [{ii.}] fix a fiber preserving map $x:T^{\ast}Q\rightarrow T^{\ast}Q$
with image inside $W$ and such that the formula $\left(\ref{cif}\right)$: 
\begin{description}
\item [{(a)}] defines a smooth application on all of $T^{\ast}Q$,
\item [{(b)}] satisfies $\left(\ref{lim}\right)$;
\end{description}
\item [{iii.}] define $z_{d}:T^{\ast}Q\rightarrow T^{\ast}Q$ by the formula
$\left(\ref{cif}\right)$; 
\end{description}
then we have a solution of $\left(\ref{szd}\right)$. Moreover, every
solution of $\left(\ref{szd}\right)$ can be constructed in this way.\end{proposition}

For instance, we can take 
\[
\mu\left(\alpha\right):=\hat{\mathfrak{b}}\left(P\left(\alpha\right),P\left(\alpha\right)\right)\ \ \ \text{and\ \ \ }x\left(\alpha\right):=-\mu\left(\alpha\right)\,\xi\left(\pi\left(\alpha\right)\right),
\]
being $\xi:Q\rightarrow T^{\ast}Q$ any $1$-form on $Q$ with image
inside $W$ (this is possible, since $W$ is a linear subbundle).

\bigskip{}

Summing up, we have a new method for solving the problem \textbf{P}.
Assume that a connection was already chosen on $T^{\ast}Q$.

\begin{definition} \label{ches}Given an underactuated system $\left(H,\mathcal{W}\right)$,
with $H=\mathfrak{H}+h\circ\pi$ simple and $\mathcal{W}$ defined
by a subbundle $W$ {[}see Eq. $\left(\ref{ww}\right)${]}, and given
a critical point $\alpha^{\bullet}\in T^{\ast}Q$ of $X_{H}$, the
\textbf{simple CH method} consists in:
\begin{itemize}
\item finding a solution $\hat{H}=\hat{\mathfrak{H}}+\hat{h}\circ\pi$ of
Eqs. $\left(\ref{kkechi}\right)$ and $\left(\ref{pechi}\right)$,
with $\hat{H}$ positive definite w.r.t. $\alpha^{\bullet}$; 
\item fixing a tensor $\mathfrak{Z}_{g}:T^{\ast}Q\times T^{\ast}Q\times T^{\ast}Q\rightarrow\mathbb{R}$
through the steps $1$ to $3$ above;
\item fixing a fiber preserving map $z_{d}:T^{\ast}Q\rightarrow T^{\ast}Q$
through the steps $i$ to $iii$ above;
\item and defining a fiber preserving map $y:T^{\ast}Q\rightarrow T^{\ast}Q$
as in $\left(\ref{y}\right)$ and $Y\in\mathfrak{X}\left(T^{\ast}Q\right)$
as the vertical lift of $y$.
\end{itemize}
\end{definition}

Of course, the simple CH method is Lyapunov based. In particular,
$\hat{H}$ and $Y$ satisfy $\left(\ref{des0}\right)$. And it is
easy to show that the simple CH method is included (in the sense of
Definition \ref{ied}) in the (general) CH method (see Definition
\ref{chm}).

\bigskip{}

Coming back to the new matching conditions $\left(\ref{kkechi}\right)$
and $\left(\ref{pechi}\right)$ {[}with local versions $\left(\ref{kec}\right)$
and $\left(\ref{pec}\right)${]}, the improvement or simplification
accomplished by Chang {[}w.r.t. to the matching conditions $\left(\ref{kei}\right)$
and $\left(\ref{pei}\right)${]} is two-fold. On the one hand, the
three unknown $\hat{\mathfrak{H}}$, $\hat{h}$ and $\mathfrak{Z}_{g}$
in $\left(\ref{kei}\right)$ and $\left(\ref{pei}\right)$ have been
\emph{decoupled}.\footnote{More precisely, we can firstly find $\hat{\mathfrak{H}}$ by solving
$\left(\ref{kkechi}\right)$, then, for such a solution $\hat{\mathfrak{H}}$,
we can find $\hat{h}$ by solving $\left(\ref{pechi}\right)$, and
finally, using $\hat{\mathfrak{H}}$ and $\hat{h}$, we can construct
$\mathfrak{Z}_{g}$ by following the steps $1$ to $3$ listed above.} On the other hand, the Eq. $\left(\ref{ke0}\right)$ {[}the local
version of $\left(\ref{kei}\right)${]} has been replaced by the Eq.
$\left(\ref{kec}\right)$, a too much simple set of equations. It
is simpler not only because of the form, but also because of the number
of equations that contains. In fact, it can be shown that the number
of equations in $\left(\ref{ke0}\right)$ and $\left(\ref{kec}\right)$
are, respectively,
\begin{equation}
\frac{n\,\left(n+1\right)\,\left(n-m\right)}{2}\label{trn}
\end{equation}
 and 
\begin{equation}
\frac{\left(n-m+2\right)\,\left(n-m+1\right)\,\left(n-m\right)}{6},\label{chn}
\end{equation}
being $n:=\dim Q$ and $m$ the rank of $W$.

\begin{remark} Regarding the last improvement, it was shown in Reference
\cite{cra} that, even for the traditional IDA-PBC method (where the
unknowns $\mathfrak{Z}_{g}^{ijk}$ adopt a particular form), the number
of equations in $\left(\ref{ke0}\right)$ can also be reduced from
$\left(\ref{trn}\right)$ to $\left(\ref{chn}\right)$ by using the
freedom one has in choosing each function $\mathfrak{Z}_{g}^{ijk}$.
Thus, the main contribution of Chang in that respect, perhaps, was
not to reduce the number of equations, but to give a precise, simple
and useful prescription to do that. \end{remark}

Equations $\left(\ref{kec}\right)$ and $\left(\ref{pec}\right)$
{[}and consequently $\left(\ref{kkechi}\right)$ and $\left(\ref{pechi}\right)${]}
were independently obtained in \cite{gym} {[}see Equations $\left(67\right)$
and $\left(68\right)$ of \cite{gym} for $\mathbb{V}^{ij}=\hat{\mathbb{H}}^{ij}$
and $v=\hat{h}${]}, almost simultaneously with the paper of Chang
\cite{chang1}, in the context of the so-called\emph{ Lyapunov constraint
based (LCB) method} for underactuated systems with only one actuator.
We will see in the last section of this paper that the same equations
are obtained (in the mentioned context) for an arbitrary number of
actuators.

\section{The LCB method}

In this section we extend the Lyapunov constraint based method for
the stabilization of underactuated systems, originally presented in
\cite{hocs} (and then further developed in \cite{gym}) for systems
with one degree of actuation, to systems with an arbitrary number
of actuators. To do that, we firstly recall, within the Hamiltonian
framework, the idea of controlling underactuated mechanical systems
by imposing kinematic constraints \cite{marle,marle2,shiriaev} (see
also \cite{g,perez1,perez2} for further examples), and the deep relationship
between constrained and closed-loop mechanical systems. It is worth
mentioning that we shall focus on a Hamiltonian formulation of the
method, although a Lagrangian one is equally possible.

\subsection{Second order constraints and closed-loop systems}

Following \cite{hocs}, a \emph{second order constrained system} (SOCS)
on $Q$ is a triple $\left(H,\mathcal{P},\mathcal{W}\right)$ where
\begin{enumerate}
\item $H:T^{\ast}Q\rightarrow\mathbb{R}$ is a smooth function defining
an (unconstrained) Hamiltonian system,
\item $\mathcal{P}\subset TT^{\ast}Q$ is a submanifold defining the \emph{second
order} \emph{kinetic constraints}\footnote{Note that, in canonical coordinates $(q,p)$, the submanifold $\mathcal{P}$
is defined, among other things, by restrictions on $\dot{q}$ and
$\dot{p}$. In the Lagrangian formalism, applying the Legendre transformation,
this gives rise to restrictions on $\dot{q}$ and $\ddot{q}$. This
is why we talk about \emph{second} order constraints.} imposed on the system, and
\item $\mathcal{W}$ is a vertical subbundle of the tangent bundle $TT^{\ast}Q$
defining the \emph{subspace of constraint forces}. 
\end{enumerate}
In this paper, by a trajectory\footnote{These curves were called \emph{type III trajectories} in \cite{hocs},
where another types of trajectories were also considered.} of $\left(H,\mathcal{P},\mathcal{W}\right)$ we mean an integral
curve $\Gamma:I\subset\mathbb{R}\rightarrow T^{\ast}Q$ of a vector
field $X\in\mathfrak{X}(T^{\ast}Q)$ that satisfies
\begin{equation}
X\subset\mathcal{P}\quad\text{and}\quad X-X_{H}\subset\mathcal{W}.\label{socstray}
\end{equation}
Of course, any trajectory must satisfy $\Gamma^{\prime}(t)\in\mathcal{P}$,
for all $t\in I$. As in the previous section, $X_{H}\in\mathfrak{X}(T^{\ast}Q)$
is the Hamiltonian vector field of $H$ w.r.t. the canonical symplectic
form of $T^{\ast}Q$. The vector field $Y:=X-X_{H}$ is called the
\emph{constraint force }related to $X$.

\begin{remark} \label{rem:Xsol}Note that $X$ is a solution of \eqref{socstray}
if and only if $Y=X-X_{H}$ is a solution of 
\begin{equation}
X_{H}+Y\subset\mathcal{P}\quad\text{and}\quad Y\subset\mathcal{W}.\label{socstray2}
\end{equation}

\end{remark}

On the other hand, a \emph{closed-loop mechanical system} (CLMS) is
defined in \cite{hocsclosed} as a dynamical system on $Q$ given
by
\begin{enumerate}
\item a smooth function $H:T^{\ast}Q\rightarrow\mathbb{R}$, describing
a (non actuated) Hamiltonian system,
\item a vertical subbundle $\mathcal{W}\subset TT^{\ast}Q$, representing
the \emph{actuation subspace} and defining, together with $H$, the
\emph{underactuated system }$\left(H,\mathcal{W}\right)$, and
\item a vector field $Y$ on $T^{\ast}Q$ such that $Y\subset\mathcal{W}$:
the \emph{control law}. 
\end{enumerate}
We will denote such a system by $(H,Y)_{\mathcal{W}}$. By a trajectory
of $(H,Y)_{\mathcal{W}}$ we mean an integral curve of the vector
field $X_{H}+Y$.

\bigskip{}

Now, let us see how a SOCS gives rise to a CLMS. Let us suppose that
a triple $\left(H,\mathcal{P},\mathcal{W}\right)$ defines a SOCS
that admits a solution $X$ of \eqref{socstray} (which is unique,
for instance, in the case of \emph{normal} SOCSs -see \cite{hocsclosed}-).
Because of Remark \ref{rem:Xsol}, this is the same as saying that
it admits a solution $Y$ of \eqref{socstray2}. Then, from the SOCS
$\left(H,\mathcal{P},\mathcal{W}\right)$ we can define the CLMS $(H,Y)_{\mathcal{W}}$,
with $Y:=X-X_{H}$. Note that
\begin{itemize}
\item both systems have the same trajectories: the integral curves of the
vector field $X=X_{H}+Y$;
\item the role of $Y$ is two fold: a constraint force for $\left(H,\mathcal{P},\mathcal{W}\right)$
and a control law for $(H,Y)_{\mathcal{W}}$. 
\end{itemize}
This construction tells us that, in order to design a control strategy
for controlling a given underactuated system $\left(H,\mathcal{W}\right)$,
we can ``think of constraints,\textquotedblright \ i.e. we can think
of the possible constraints $\mathcal{P}$ that give rise to the desirable
behavior, and then obtain the control law as the related constraint
force $Y\subset\mathcal{W}$.

It was shown in \cite{hocsclosed} that every CLMS can be constructed
from a SOCS as we did above, i.e. every control law may be seen as
the constraint force of a given set of second order constraints. This
result reveals a deep connection between closed-loop and constrained
mechanical systems and, from the point of view of the applications
to automatic control, the result says that, in order to synthesize
a state feedback for a given underactuated system, we \textbf{always}\ (i.e.
without loss of generality) can ``think of constraints.\textquotedblright{}

\subsection{(Asymptotic) stability and related constraints}

Let us consider a dynamical system on a manifold $P$ defined by a
vector field $X\in\mathfrak{X}(P)$. Given a critical point $\alpha^{\bullet}\in P$
of $X$, let $\hat{H}:P\rightarrow\mathbb{R}$ be a Lyapunov function
for $X$ and $\alpha^{\bullet}$ (recall Definition \ref{lyap}).
Note that, given a trajectory $\Gamma:I\subset\mathbb{R}\rightarrow P$
of $X$, condition \textbf{L2 }of Definition \ref{lyap} implies that
\begin{equation}
\left\langle d\hat{H}(\Gamma(t)),\Gamma^{\prime}(t)\right\rangle =-\mu(\Gamma(t)),\:\:\:\forall t\in I,\label{vinc2}
\end{equation}
where $\mu:P\rightarrow\mathbb{R}$ is the non-negative function given
by
\[
\mu\left(\alpha\right):=-\left\langle d\hat{H}(\alpha),X(\alpha)\right\rangle ,\ \ \ \forall\alpha\in P.
\]

\begin{remark} Observe that $\mu^{-1}\left(0\right)$ is the \emph{La'Salle
surface} related to $\hat{H}$ (see \cite{khalil}), and $\alpha^{\bullet}\in\mu^{-1}\left(0\right)$.
\end{remark}

That is, condition \textbf{L2} may be interpreted as a kinematic constraint
on the system. Hence, roughly speaking, if we want to stabilize a
dynamical system, we can think of imposing a constraint of the form
\eqref{vinc2}, for appropriate non-negative functions $\hat{H}$
and $\mu$. We shall call it \textbf{Lyapunov constraint}. Of course,
depending on the conditions we impose on $\hat{H}$ and $\mu$, we
shall have different stability properties. For instance if, besides
condition \textbf{L1} for $\hat{H}$, we ask $\mu$ to be such that
the singleton $\left\{ \alpha^{\bullet}\right\} $ is the bigger invariant
subset of $\mu^{-1}\left(0\right)$, the \emph{La'Salle invariance
principle} would ensure (local) asymptotic stability for $\alpha^{\bullet}$.
This is true, for example, if we assume that property \textbf{L1}
also holds for $\mu$, what would imply that
\[
\left\langle d\hat{H}(\alpha),X(\alpha)\right\rangle <0\ \ \ \text{for all\ \ \ }\alpha\neq\alpha^{\bullet}.
\]
If in addition we ask $\hat{H}$ to be a proper function (and $P$
to be connected), then global asymptotic stability for $\alpha^{\bullet}$
would be ensured. (For a proof of these results, see \cite{khalil}
again.)

\bigskip{}

Now, let us focus our attention on Hamiltonian systems. Take $P=T^{\ast}Q$
for some $Q$, fix a smooth function $H:T^{\ast}Q\rightarrow\mathbb{R}$
and consider the Hamiltonian system on $Q$ defined by $H$. Given
a point $\alpha^{\bullet}\in T^{\ast}Q$ and non-negative functions
$\hat{H},\mu:T^{\ast}Q\rightarrow\mathbb{R}$, let us impose the constraint
\eqref{vinc2} on this system. In other words, let us define the submanifold
\[
\mathcal{P}:=\left\{ V\in TT^{\ast}Q:\left\langle d\hat{H}(\tau_{T^{\ast}Q}\left(V\right)),V\right\rangle =-\mu(\tau_{T^{\ast}Q}\left(V\right))\right\} ,
\]
and impose the constraint $\Gamma^{\prime}(t)\in\mathcal{P}$ on the
trajectories.

\begin{remark} Notice that, if $\left(U,\varphi\right)$ is a coordinate
chart of $Q$, in terms of the induced chart on $TT^{\ast}Q$ (see
Eq. \eqref{icor2}) the submanifold $\mathcal{P}$ is locally given
by the equation 
\[
\frac{\partial\hat{H}}{\partial q^{i}}(\mathbf{q},\mathbf{p})\,\dot{q}^{i}+\frac{\partial\hat{H}}{\partial p_{i}}(\mathbf{q},\mathbf{p})\,\dot{p}_{i}=-\mu(\mathbf{q},\mathbf{p}).
\]

\end{remark}

Suppose that we want to implement this constraint by exerting forces
lying inside a vertical subbundle $\mathcal{W}\subset TT^{\ast}Q$.
All that defines the SOCS $(H,\mathcal{P},\mathcal{W})$. Assume that
this SOCS admits a solution $X$ of \eqref{socstray}, or equivalently,
admits a solution $Y$ of \eqref{socstray2}. (In other words, assume
that the Lyapunov constraint $\mathcal{P}$ can be implemented by
a constraint force $Y\subset\mathcal{W}$.) This is the same as saying
that there exists $Y\in\mathfrak{X}\left(T^{\ast}Q\right)$ such that
\begin{equation}
\left\langle d\hat{H}(\alpha),X_{H}(\alpha)+Y(\alpha)\right\rangle =-\mu(\alpha)\ \ \ \text{and\ \ \ }Y(\alpha)\in\mathcal{W},\label{vincf0}
\end{equation}
or equivalently
\begin{equation}
\mathfrak{i}_{X_{H}+Y}\ d\hat{H}=-\mu\ \ \ \text{and\ \ \ }Y\subset\mathcal{W}.\label{vincf1}
\end{equation}
Since $\mu$ is non-negative, then $\left\langle d\hat{H}(\alpha),X_{H}(\alpha)+Y(\alpha)\right\rangle \leq0$,
i.e. $\hat{H}$ satisfies \textbf{L2}. In addition, if $X_{H}(\alpha^{\bullet})+Y(\alpha^{\bullet})=0$
and $\hat{H}$ satisfies \textbf{L1} for $\alpha^{\bullet}$, then
$\hat{H}$ is a Lyapunov function for $X_{H}+Y$ and $\alpha^{\bullet}$,
and consequently the underactuated system $\left(H,\mathcal{W}\right)$
can be stabilized at $\alpha^{\bullet}$ by the control law $Y$.
Of course, if stronger conditions are imposed on $\hat{H}$ and $\mu$
(as discussed at the beginning of this section), stronger stability
properties can be ensured for the system defined by $X_{H}+Y$.

\begin{remark} In this way, as we have seen in the previous section,
we are constructing the CLMS $(H,Y)_{\mathcal{W}}$ from the SOCS
$(H,\mathcal{P},\mathcal{W})$. \end{remark}

In conclusion, if a solution $Y$ exists for Equation \eqref{vincf0},
for some functions $\hat{H}$ and $\mu$, different assertions about
the stabilizability around $\alpha^{\bullet}$ of the underactuated
system $\left(H,\mathcal{W}\right)$ can be made, depending on the
properties of $\hat{H}$ and $\mu$.

\begin{remark} Also, if a solution $Y$ of \eqref{vincf0} exists
along an open subset $\pi^{-1}\left(U\right)=T^{\ast}U\subset T^{\ast}Q$
containing $\alpha^{\bullet}$, namely a \emph{local solution} of
\eqref{vincf0}, the same assertions can be made, just replacing $Q$
by $U$. \end{remark}

\subsection{A \emph{maximal} stabilization method}

The discussion in the previous section drives us to another method
for (asymptotic) stabilization of non-linear underactuated mechanical
systems.

\begin{definition} \label{dlm}Given an underactuated system $\left(H,\mathcal{W}\right)$
on $Q$ and a critical point $\alpha^{\bullet}\in T^{\ast}Q$ of $X_{H}$,
the \textbf{Lyapunov constraint based (LCB) method} consists in finding
two functions $\hat{H},\mu:T^{\ast}Q\rightarrow\mathbb{R}$ and a
vector field $Y\in\mathfrak{X}\left(T^{\ast}Q\right)$ such that $\hat{H}$
is positive definite w.r.t. $\alpha^{\bullet}$, $\mu$ is non-negative,
$Y(\alpha^{\bullet})=0$ and Eq. $\left(\ref{vincf0}\right)$ is solved.
\end{definition}

Note that the method is a Lyapunov based stabilization method, and
it is essentially defined by Eq. $\left(\ref{vincf0}\right)$. So,
we can identify the method with this equation. Let us write it in
other terms. Since $\mathcal{W}$ is a vertical subbundle and $Y$
is a vertical vector field, we can write $\mathcal{W}_{\alpha}=\operatorname{vlift}_{\alpha}\left(W_{\alpha}\right)$
and $Y\left(\alpha\right)=\operatorname{vlift}_{\alpha}\left(y\left(\alpha\right)\right)$,
for a unique subspace $W_{\alpha}\subset T_{\pi(\alpha)}^{\ast}Q$
and a unique fiber preserving map $y:T^{\ast}Q\rightarrow T^{\ast}Q$.
Under this notation,
\[
\left\langle d\hat{H}(\alpha),Y(\alpha)\right\rangle =\left\langle y(\alpha),\mathbb{F}\hat{H}(\alpha)\right\rangle 
\]
and conditions $Y(\alpha^{\bullet})=0$ and $Y(\alpha)\in\mathcal{W}_{\alpha}$
translate to $y(\alpha^{\bullet})=0$ and $y(\alpha)\in W_{\alpha}$,
respectively. On the other hand, 
\[
\left\langle d\hat{H}(\alpha),X_{H}(\alpha)\right\rangle =\{\hat{H},H\}(\alpha),
\]
being $\{\hat{H},H\}$ the canonical Poisson bracket between $\hat{H}$
and $H$. Combining all that, we can write $\left(\ref{vincf0}\right)$
as
\begin{equation}
\left\langle y(\alpha),\mathbb{F}\hat{H}(\alpha)\right\rangle =-\mu(\alpha)-\left\{ \hat{H},H\right\} (\alpha)\ \ \ \text{and\ \ \ }y(\alpha)\in W.\label{vincf}
\end{equation}

\bigskip{}

To conclude the section, let us mention the important fact that any
stabilization method for $\left(H,\mathcal{W}\right)$ which gives
rise to a control law $Y\subset\mathcal{W}$ and a Lyapunov function
$\hat{H}$ for the related closed-loop system $X_{H}+Y$ (and some
critical point of $X_{H}$), as every version of the energy shaping
method does, can be reduced to the LCB method, i.e. to solve Eq. \eqref{vincf0}
{[}or equivalently, to solve Eq. \eqref{vincf}{]}. More precisely,

\begin{theorem} \label{inc}Let $\left(H,\mathcal{W}\right)$ be
an underactuated system and $\alpha^{\bullet}\in T^{\ast}Q$ a critical
point of $X_{H}$. If we are given a vector field $Y\subset\mathcal{W}$
and a Lyapunov function $\hat{H}$ for $X_{H}+Y$ and $\alpha^{\bullet}$,
then $\hat{H}$ is positive definite w.r.t. $\alpha^{\bullet}$, $\mu:=-\mathfrak{i}_{X_{H}+Y}\ d\hat{H}$
is non-negative and $Y\left(\alpha^{\bullet}\right)=0$. In particular,
$Y$ is given by the LCB method. \end{theorem}

\emph{Proof.} Given a vector field $Y\subset\mathcal{W}$, if $\hat{H}$
is a Lyapunov function for $X:=X_{H}+Y$ and $\alpha^{\bullet}$,
the theorem easily follows from the fact that $\alpha^{\bullet}$
must be critical for $X$ {[}from which $Y\left(\alpha^{\bullet}\right)=0${]},
the item \textbf{L1} and the combination of \textbf{L2} and the Eq.
\eqref{vincf1}.\ \ \ $\square$

\bigskip{}

In terms of Definitions \ref{smd} and \ref{ied}, this theorem says
that the LCB method includes all the Lyapunov based stabilization
methods, i.e. it is \textbf{\emph{maximal}} among such methods. In
other words, the LCB method is the most general method among the Lyapunov
based stabilization methods (see Remark \ref{generality}). In particular,
any version of the CH method (see Definition \ref{chm}) is included
in the LCB method.

\section{The LCB method for simple functions}

In Ref. \cite{gym}, a deep study of Eq. \eqref{vincf} has been done
for underactuated systems with only one actuator. In this section
we shall extend such a study to an arbitrary number of actuators.
More precisely, given an underactuated system $\left(H,\mathcal{W}\right)$
and non-negative functions $\hat{H}$ and $\mu$, we shall study under
which conditions there exists a fiber preserving map $y$ solving
\eqref{vincf} (thinking of $\hat{H}$ and $\mu$ as data of \eqref{vincf},
instead of unknowns). We shall focus in the case in which $H$ and
$\hat{H}$ are simple functions. In this case, we show that the mentioned
existence problem is governed by a set of PDEs for $\hat{H}$, which
define what we have called the \textbf{\emph{simple}}\emph{ LCB method}
in Ref. \cite{gym}. We shall see that these equations are exactly
the matching conditions $\left(\ref{kkechi}\right)$ and $\left(\ref{pechi}\right)$
obtained by Chang in \cite{chang3,chang4} (related to the simple
CH method of Definition \ref{ches}). Finally, we show that the simple
LCB and the simple CH method are equivalent stabilization methods.

\subsection{The \emph{kinetic }and\emph{ potential }equations}

Following the same notation as in Section \ref{kpmc}, consider two
simple Hamiltonian functions $H=\mathfrak{H}+h\circ\pi$ and $\hat{H}=\hat{\mathfrak{H}}+\hat{h}\circ\pi$
on the cotangent bundle $T^{\ast}Q$. Then, given a torsion-free connection
on $T^{\ast}Q$, the canonical Poisson bracket between $H$ and $\hat{H}$
can be written, for all $\alpha\in T^{\ast}Q$,
\[
\left\{ \hat{H},H\right\} \left(\alpha\right)=\left\langle \mathbb{B}\hat{\mathfrak{H}}\left(\alpha\right),\mathbb{F}\mathfrak{H}\left(\alpha\right)\right\rangle -\left\langle \mathbb{B}\mathfrak{H}\left(\alpha\right),\mathbb{F}\hat{\mathfrak{H}}\left(\alpha\right)\right\rangle +\left\langle d\hat{h}\left(q\right),\mathbb{F}\mathfrak{H}\left(\alpha\right)\right\rangle -\left\langle dh\left(q\right),\mathbb{F}\hat{\mathfrak{H}}\left(\alpha\right)\right\rangle ,
\]
with $q=\pi\left(\alpha\right)$. This is a direct consequence of
Eqs. $\left(\ref{basic}\right)$ and $\left(\ref{fiberbasebracket}\right)$.
In a local chart $\left(U,\varphi\right)$, using $\left(\ref{id1}\right)$
and $\left(\ref{id2}\right)$ (and the fact that our connection is
torsion-free), 
\begin{equation}
\left\{ \hat{H},H\right\} (\mathbf{q},\mathbf{p})=\frac{1}{2}\left(\frac{\partial\hat{\mathbb{H}}^{ij}(\mathbf{q})}{\partial q^{k}}\,\mathbb{H}^{kl}(\mathbf{q})-\frac{\partial\mathbb{H}^{ij}(\mathbf{q})}{\partial q^{k}}\,\hat{\mathbb{H}}^{kl}(\mathbf{q})\right)\,p_{i}p_{j}p_{l}+\left(\frac{\partial\hat{h}(\mathbf{q})}{\partial q^{k}}\,\mathbb{H}^{kl}(\mathbf{q})-\frac{\partial h(\mathbf{q})}{\partial q^{k}}\,\hat{\mathbb{H}}^{kl}(\mathbf{q})\right)\,p_{l}.\label{cpb}
\end{equation}
Looking at above expression, the next result easily follows.

\begin{lemma} If $H$ and $\hat{H}$ are simple functions on $T^{\ast}Q$,
then $\left\{ \hat{H},H\right\} $ is an odd function (when restricted
to each fiber of $T^{\ast}Q$), i.e. 
\begin{equation}
\left\{ \hat{H},H\right\} \left(-\alpha\right)=-\left\{ \hat{H},H\right\} \left(\alpha\right)\ \ \ \ \ \text{for\ all\ \ \ }\alpha\in T^{\ast}Q.\label{odd}
\end{equation}

\end{lemma}

Now, fix an underactuated system $\left(H,\mathcal{W}\right)$ with
$H$ simple and $\mathcal{W}$ given by a subbundle $W\subset T^{\ast}Q$
{[}see $\left(\ref{ww}\right)${]}.

\begin{proposition} \label{prop1}Consider a simple function $\hat{H}:T^{\ast}Q\rightarrow\mathbb{R}$
and a non-negative function $\mu:T^{\ast}Q\rightarrow\mathbb{R}$.
If there exists a solution $y:T^{\ast}Q\rightarrow T^{\ast}Q$ of
Eq. \eqref{vincf} for $\hat{H}$ and $\mu$, then 
\begin{equation}
\left\{ \hat{H},H\right\} \left(\sigma\right)=0,\ \ \ \ \forall\sigma\in\hat{W},\label{t0}
\end{equation}
where $\hat{W}:=\left(\mathbb{F}\hat{\mathfrak{H}}\left(W\right)\right)^{0}=\mathbb{F}\hat{\mathfrak{H}}^{-1}\left(W^{0}\right)$.
\end{proposition}

\emph{Proof.} By hypothesis {[}see Eq. \eqref{vincf}{]}
\[
\left\langle \alpha,\mathbb{F}\hat{\mathfrak{H}}(y(\alpha))\right\rangle =-\mu(\alpha)-\{\hat{H},H\}(\alpha)\ \ \ \text{and\ \ \ }y(\alpha)\in W.
\]
It is clear that $\left\langle \sigma,\mathbb{F}\hat{\mathfrak{H}}(y(\sigma))\right\rangle =0$
for all $\sigma\in\hat{W}$, and consequently
\[
\mu(\sigma)+\{\hat{H},H\}(\sigma)=0,\ \ \ \ \forall\sigma\in\hat{W}.
\]
Suppose that for some $\sigma_{0}\in\hat{W}$ we have that $\{\hat{H},H\}(\sigma_{0})\neq0$.
If $\{\hat{H},H\}(\sigma_{0})>0$, then $\mu(\sigma_{0})<0$. But
this is not possible, so $\{\hat{H},H\}(\sigma_{0})<0$. According
to Eq. \eqref{odd}, 
\[
\{\hat{H},H\}(-\sigma_{0})=-\{\hat{H},H\}(\sigma_{0})>0.
\]
Since $-\sigma_{0}\in\hat{W}$, it follows that $\mu(-\sigma_{0})+\{\hat{H},H\}(-\sigma_{0})=0$,
what implies $\mu(-\sigma_{0})<0$. As a consequence, $\{\hat{H},H\}(\sigma)=0$
for all $\sigma\in\hat{W}$.\ \ \ $\square$

\bigskip{}

We can write $\left(\ref{t0}\right)$ in an equivalent way.

\begin{lemma} \label{equivs}If $H$ and $\hat{H}$ are simple functions
on $T^{\ast}Q$ and $V\subset T^{\ast}Q$ is a linear subbundle, the
following conditions are equivalent.
\begin{enumerate}
\item $\left\{ \hat{H},H\right\} \left(\sigma\right)=0$, $\forall\sigma\in V$.
\item Given a connection on $T^{\ast}Q$, for all $q\in Q$ and $\sigma\in V_{q}$,
the Eqs. $\left(\ref{kkechi}\right)$ and $\left(\ref{pechi}\right)$
hold, with $\hat{W}$ replaced by $V$.
\item Given a local chart $\left(U,\varphi\right)$, for all $q\in U$ and
$\mathbf{p}\in\left(\varphi_{q}^{\ast}\right)^{-1}\left(V_{q}\right)$,
the Eqs. $\left(\ref{kec}\right)$ and $\left(\ref{pec}\right)$ hold. 
\end{enumerate}
\end{lemma}

\emph{Proof.} The equivalence between $2$ and $3$ is given by the
Eqs. $\left(\ref{id1}\right)$ and $\left(\ref{id2}\right)$. We only
need to prove that $1$ implies $3$ (the converse is immediate).
If $\{\hat{H},H\}\left(\sigma\right)=0$, $\forall\sigma\in V$, it
follows from Equation $\left(\ref{cpb}\right)$ that, in a local chart
$\left(U,\varphi\right)$,
\begin{equation}
\frac{1}{2}\left(\frac{\partial\hat{\mathbb{H}}^{ij}(\mathbf{q})}{\partial q^{k}}\,\mathbb{H}^{kl}(\mathbf{q})-\frac{\partial\mathbb{H}^{ij}(\mathbf{q})}{\partial q^{k}}\,\hat{\mathbb{H}}^{kl}(\mathbf{q})\right)\,p_{i}p_{j}p_{l}+\left(\frac{\partial\hat{h}(\mathbf{q})}{\partial q^{k}}\,\mathbb{H}^{kl}(\mathbf{q})-\frac{\partial h(\mathbf{q})}{\partial q^{k}}\,\hat{\mathbb{H}}^{kl}(\mathbf{q})\right)\,p_{l}=0,\label{twoterms}
\end{equation}
for all $q\in U$ and $\mathbf{p}\in\left(\varphi_{q}^{\ast}\right)^{-1}\left(V_{q}\right)$.
We must show that each term vanishes. Fixing $\mathbf{\widehat{p}}\in\left(\varphi_{q}^{\ast}\right)^{-1}\left(V_{q}\right)$,
define 
\[
A:=\left(\frac{\partial\hat{\mathbb{H}}^{ij}(\mathbf{q})}{\partial q^{k}}\,\mathbb{H}^{kl}(\mathbf{q})-\frac{\partial\mathbb{H}^{ij}(\mathbf{q})}{\partial q^{k}}\,\hat{\mathbb{H}}^{kl}(\mathbf{q})\right)\,\,\widehat{p}_{i}\widehat{p}_{j}\widehat{p}_{l}\ \ \ \text{and\ \ \ }B:=\left(\frac{\partial\hat{h}(\mathbf{q})}{\partial q^{k}}\,\mathbb{H}^{kl}(\mathbf{q})-\frac{\partial h(\mathbf{q})}{\partial q^{k}}\,\hat{\mathbb{H}}^{kl}(\mathbf{q})\right)\,\widehat{p}_{l}.
\]
Then, if we replace $\mathbf{p}$ by $\lambda\,\mathbf{\widehat{p}}$
in Eq. \eqref{twoterms}, with $\lambda\in\mathbb{R}$, we have that
$A\,\lambda^{3}/2+B\,\lambda=0$ for all $\lambda\in\mathbb{R}$.
But this is possible only if $A=B=0$. \ \ \ $\square$

\bigskip{}

Combining last lemma and Proposition \ref{prop1}, we have the following
result.

\begin{proposition} \label{prop2}Consider a simple function $\hat{H}:T^{\ast}Q\rightarrow\mathbb{R}$
and a non-negative function $\mu:T^{\ast}Q\rightarrow\mathbb{R}$.
If there exists a solution $y:T^{\ast}Q\rightarrow T^{\ast}Q$ of
Eq. \eqref{vincf} for $\hat{H}$ and $\mu$, then, for every connection
on $T^{\ast}Q$,
\begin{align}
\left\langle \mathbb{B}\hat{\mathfrak{H}}\left(\sigma\right),\mathbb{F}\mathfrak{H}\left(\sigma\right)\right\rangle -\left\langle \mathbb{B}\mathfrak{H}\left(\sigma\right),\mathbb{F}\hat{\mathfrak{H}}\left(\sigma\right)\right\rangle  & =0,\label{lcbm1}\\
\left\langle d\hat{h}\left(\pi\left(\sigma\right)\right),\mathbb{F}\mathfrak{H}\left(\sigma\right)\right\rangle -\left\langle dh\left(\pi\left(\sigma\right)\right),\mathbb{F}\hat{\mathfrak{H}}\left(\sigma\right)\right\rangle  & =0,\ \ \ \forall\sigma\in\hat{W}.\label{lcbm}
\end{align}

\end{proposition}

Notice that Equations $\left(\ref{lcbm1}\right)$ and $\left(\ref{lcbm}\right)$
are exactly the \emph{simple} matching conditions $\left(\ref{kkechi}\right)$
and $\left(\ref{pechi}\right)$ derived by Chang. It is quite surprising
for us that these two methods, the LCB and the simple CH methods,
give rise to the same set of equations, in spite of they arise from
very different ideas: ``feedback equivalence'' and ``controlling
by imposing kinematic constraints.''

\begin{remark} This means, using again Lemma \ref{equivs}, that
the simple matching conditions can be written as in $\left(\ref{t0}\right)$,
or equivalently 
\[
\left\{ \hat{H},H\right\} \circ\mathbb{F}\hat{H}^{-1}\left(v\right)=0,\ \ \ \ \forall v\in W^{0},
\]
which is a coordinate-free and connection-free equation. \end{remark}

Propositions \ref{prop1} and \ref{prop2} say that Equation $\left(\ref{t0}\right)$,
or equivalently Equations $\left(\ref{lcbm1}\right)$ and $\left(\ref{lcbm}\right)$,
which have a function $\hat{H}=\hat{\mathfrak{H}}+\hat{h}\circ\pi$
as unknown, define necessary conditions for the existence of a solution
of \eqref{vincf0} for some non-negative function $\mu$. Let us see
that they also define sufficient conditions. We saw in Section \ref{cmc}
that if a simple function $\hat{H}$ satisfies $\left(\ref{lcbm1}\right)$
and $\left(\ref{lcbm}\right)$, then there exists a vector field $Y\subset\mathcal{W}$
satisfying $\left(\ref{des0}\right)$. (In fact, $Y$ can be constructed
by following the procedure described in Definition \ref{ches}). Thus,
$Y$ satisfies Eq. \eqref{vincf0} for $\hat{H}$ and some non-negative
function $\mu$. Concluding,

\begin{theorem} \label{ssi}Given a simple function $\hat{H}:T^{\ast}Q\rightarrow\mathbb{R}$,
there exists a solution $y$ of Eq. \eqref{vincf}, for $\hat{H}$
and some non-negative function $\mu$, if and only if $\left(\ref{t0}\right)$
is fulfilled for $\hat{H}$. \end{theorem}

\begin{remark} In contrast to the situation considered in Ref. \cite{gym},
since here the subbundle $W$ is not $1$-dimensional in general,
we can not ensure the uniqueness of solutions $y$ of Eq. \eqref{vincf}
(for $\hat{H}$ and $\mu$ fixed). It is shown in \cite{gym} that,
if a global generator $\xi:Q\rightarrow T^{\ast}Q$ of the $1$-dimensional
subbundle $W$ is given, the solution is necessarily
\[
y\left(\alpha\right):=-\frac{\mu\left(\alpha\right)+\left\{ \hat{H},H\right\} \left(\alpha\right)}{\left\langle \xi\left(\pi\left(\alpha\right)\right),\mathbb{F}\hat{H}\left(\alpha\right)\right\rangle }\ \xi\left(\pi\left(\alpha\right)\right).
\]

\end{remark}

Taking into account the maximality of the LCB method (see Theorem
\ref{inc}), the theorem above says that $\left(\ref{t0}\right)$,
and consequently the simple matching conditions, are a \textbf{\emph{key
ingredient}} that must be present (explicitly or not) in any stabilization
method that involves the construction of a simple Lyapunov function.
More precisely, if by means of some stabilization method we construct
a control law $Y$ and a simple Lyapunov function $\hat{H}$ for $X_{H}+Y$,
then $\hat{H}$ must satisfy Eq. $\left(\ref{t0}\right)$.

\subsection{The simple LCB method}

Let us see how the procedure described in Definition \ref{dlm} can
be reformulated when simple Hamiltonian functions are involved. Fix
an underactuated system $\left(H,\mathcal{W}\right)$ with $H$ simple
and $\mathcal{W}$ given by a subbundle $W\subset T^{\ast}Q$. Fix
also a critical point $\alpha^{\bullet}=\left(q^{\bullet},0\right)$
of $X_{H}$ (see Remark \ref{cps}). Suppose that we have a function
$\hat{H}$, simple and positive-definite w.r.t. $\alpha^{\bullet}$,
and a fiber preserving map $y:T^{\ast}Q\rightarrow T^{\ast}Q$ such
that $y\left(\alpha^{\bullet}\right)=0$. It is clear that $y$ can
be written
\begin{equation}
y\left(\alpha\right):=z\left(\alpha\right)-\mathbb{F}\hat{H}^{-1}\circ\left(\mathbb{F}H\circ\mathbb{B}\hat{H}\left(\alpha\right)-\mathbb{F}\hat{H}\circ\mathbb{B}H\left(\alpha\right)\right),\label{y1}
\end{equation}
being $z:T^{\ast}Q\rightarrow T^{\ast}Q$ some fiber preserving map.
Consider the orthogonal decomposition $T^{\ast}Q=W\oplus\hat{W}$
(recall Remark \ref{ort1}), given by the tensor $\hat{\mathfrak{b}}$,
and write $z\left(\alpha\right)=z_{\shortparallel}\left(\alpha\right)-z_{\perp}\left(\alpha\right)$,
with $z_{\shortparallel}\left(\alpha\right)\in W$ and $z_{\perp}\left(\alpha\right)\in\hat{W}$.

\begin{proposition}Under above assumptions and notation, the map
$y$ given by $\left(\ref{y1}\right)$ satisfies $\left(\ref{vincf}\right)$,
for some non-negative function $\mu$, if and only if $\hat{H}$ satisfies
$\left(\ref{t0}\right)$ and $z$ satisfies
\begin{equation}
\hat{\mathfrak{b}}\left(z\left(\alpha\right),\alpha\right)=-\mu\left(\alpha\right),\ \ \ \hat{\mathfrak{b}}\left(z\left(\alpha\right),\sigma\right)=\Upsilon\left(\alpha,\alpha,\sigma\right)\ \ \ \ \text{and}\ \ \ z\left(\alpha^{\bullet}\right)=0,\label{fsp}
\end{equation}
{[}where $\Upsilon:T^{\ast}Q\times T^{\ast}Q\times T^{\ast}Q\rightarrow\mathbb{R}$
is the tensor field defined in $\left(\ref{T}\right)${]} for all
$\alpha\in T^{\ast}Q$ and $\sigma\in\hat{W}$; or equivalently, $z_{\perp}$
is given by
\begin{equation}
\hat{\mathfrak{b}}\left(z_{\perp}\left(\alpha\right),\sigma\right)=\left\{ \begin{array}{rr}
0, & \ \ \ \sigma\in W,\\
-\Upsilon\left(\alpha,\alpha,\sigma\right), & \ \ \ \sigma\in\hat{W},
\end{array}\right.\label{zo}
\end{equation}
and $z_{\shortparallel}$ fulfills
\begin{equation}
\hat{\mathfrak{b}}\left(z_{\shortparallel}\left(\alpha\right),\alpha\right)=\hat{\mathfrak{b}}\left(z_{\perp}\left(\alpha\right),\alpha\right)-\mu\left(\alpha\right),\ \ z_{\shortparallel}\left(\alpha\right)\in W\ \ \ \text{and\ \ \ }z_{\shortparallel}\left(\alpha^{\bullet}\right)=0.\label{zos}
\end{equation}
\end{proposition}

\emph{Proof. }Let us first recall that, from the results of the previous
section (see Theorem \ref{ssi}), there exists a solution $y$ of
$\left(\ref{vincf}\right)$ for $\hat{H}$ (and for some non-negative
function $\mu$) if and only if $\hat{H}$ satisfies $\left(\ref{t0}\right)$,
or equivalently $\left(\ref{lcbm1}\right)$ and $\left(\ref{lcbm}\right)$.
If such a solution is given by $\left(\ref{y1}\right)$, let us see
what that means in terms of $z$. Using $\left(\ref{fiberbasebracket}\right)$,
we have that
\[
\begin{array}{l}
\left\langle \mathbb{F}\hat{H}^{-1}\circ\left(\mathbb{F}H\circ\mathbb{B}\hat{H}\left(\alpha\right)-\mathbb{F}\hat{H}\circ\mathbb{B}H\left(\alpha\right)\right),\mathbb{F}\hat{H}\left(\alpha\right)\right\rangle =\left\langle \alpha,\mathbb{F}H\circ\mathbb{B}\hat{H}\left(\alpha\right)-\mathbb{F}\hat{H}\circ\mathbb{B}H\left(\alpha\right)\right\rangle =\\
\\
=\left\langle \mathbb{B}\hat{H}\left(\alpha\right),\mathbb{F}H\left(\alpha\right)\right\rangle -\left\langle \mathbb{B}H\left(\alpha\right),\mathbb{F}\hat{H}\left(\alpha\right)\right\rangle =\left\{ \hat{H},H\right\} \left(\alpha\right).
\end{array}
\]
Then $\left\langle y\left(\alpha\right),\mathbb{F}\hat{H}\left(\alpha\right)\right\rangle =\left\langle z\left(\alpha\right),\mathbb{F}\hat{H}\left(\alpha\right)\right\rangle -\left\{ \hat{H},H\right\} \left(\alpha\right)$
and, consequently, the first part of $\left(\ref{vincf}\right)$ implies
that 
\begin{equation}
\left\langle z\left(\alpha\right),\mathbb{F}\hat{H}\left(\alpha\right)\right\rangle =-\mu\left(\alpha\right).\label{fp}
\end{equation}
Let us study the second part of $\left(\ref{vincf}\right)$, i.e.
the condition $y\subset W$. First note that, for all $\alpha,\sigma\in T^{\ast}Q$
on the same fiber,
\begin{equation}
\left\langle y\left(\alpha\right),\mathbb{F}\hat{H}\left(\sigma\right)\right\rangle =\left\langle z\left(\alpha\right),\mathbb{F}\hat{H}\left(\sigma\right)\right\rangle -\left\langle \mathbb{B}\hat{H}\left(\alpha\right),\mathbb{F}H\left(\sigma\right)\right\rangle +\left\langle \mathbb{B}H\left(\alpha\right),\mathbb{F}\hat{H}\left(\sigma\right)\right\rangle .\label{yf}
\end{equation}
To find a more convenient expression of the above equation, let us
write $H$ and $\hat{H}$ in terms of their kinetic and potential
terms, and consider the tensor $\Upsilon$ defined in $\left(\ref{T}\right)$
and the tensor $\psi:T^{\ast}Q\rightarrow\mathbb{R}$ given by 
\[
\psi\left(\alpha\right):=\left\langle d\hat{h}\left(\pi\left(\alpha\right)\right),\mathbb{F}\mathfrak{H}\left(\cdot\right)\right\rangle -\left\langle dh\left(\pi\left(\alpha\right)\right),\mathbb{F}\hat{\mathfrak{H}}\left(\cdot\right)\right\rangle .
\]
According to the second part of $\left(\ref{bq}\right)$, Eq. $\left(\ref{T}\right)$
says that
\[
\Upsilon\left(\alpha,\alpha,\sigma\right)=\left\langle \mathbb{B}\hat{\mathfrak{b}}\left(\alpha,\alpha\right),\mathbb{F}\mathfrak{H}\left(\sigma\right)\right\rangle -\left\langle \mathbb{B}\mathfrak{b}\left(\alpha,\alpha\right),\mathbb{F}\hat{\mathfrak{H}}\left(\sigma\right)\right\rangle =\left\langle \mathbb{B}\hat{\mathfrak{H}}\left(\alpha\right),\mathbb{F}\mathfrak{H}\left(\sigma\right)\right\rangle -\left\langle \mathbb{B}\mathfrak{H}\left(\alpha\right),\mathbb{F}\hat{\mathfrak{H}}\left(\sigma\right)\right\rangle ,
\]
for all $\alpha,\sigma\in T^{\ast}Q$ on the same fiber. Then, Eqs.
$\left(\ref{lcbm1}\right)$ and $\left(\ref{lcbm}\right)$ read 
\[
\Upsilon\left(\sigma,\sigma,\sigma\right)=\psi\left(\sigma\right)=0,\ \ \ \forall\sigma\in\hat{W}.
\]
On the other hand, since
\[
\begin{array}{l}
\left\langle \mathbb{B}\hat{H}\left(\alpha\right),\mathbb{F}H\left(\sigma\right)\right\rangle -\left\langle \mathbb{B}H\left(\alpha\right),\mathbb{F}\hat{H}\left(\sigma\right)\right\rangle =\left\langle \mathbb{B}\hat{\mathfrak{H}}\left(\alpha\right),\mathbb{F}\mathfrak{H}\left(\sigma\right)\right\rangle -\left\langle \mathbb{B}\mathfrak{H}\left(\alpha\right),\mathbb{F}\hat{\mathfrak{H}}\left(\sigma\right)\right\rangle +\\
\\
+\left(\left\langle d\hat{h}\left(q\right),\mathbb{F}\mathfrak{H}\left(\cdot\right)\right\rangle -\left\langle dh\left(q\right),\mathbb{F}\hat{\mathfrak{H}}\left(\cdot\right)\right\rangle \right)=\Upsilon\left(\alpha,\alpha,\sigma\right)+\psi\left(\sigma\right),
\end{array}
\]
we have from $\left(\ref{yf}\right)$ the equality
\begin{equation}
\left\langle y\left(\alpha\right),\mathbb{F}\hat{H}\left(\sigma\right)\right\rangle =\left\langle z\left(\alpha\right),\mathbb{F}\hat{H}\left(\sigma\right)\right\rangle -\Upsilon\left(\alpha,\alpha,\sigma\right)-\psi\left(\sigma\right).\label{yf1}
\end{equation}
Coming back to $\left(\ref{vincf}\right)$, the condition $y\subset W$
implies that, for all $\sigma\in\hat{W}=\mathbb{F}\hat{H}^{-1}\left(W^{0}\right)$
{[}see Eq. $\left(\ref{yf1}\right)${]}, 
\[
0=\left\langle y\left(\alpha\right),\mathbb{F}\hat{H}\left(\sigma\right)\right\rangle =\left\langle z\left(\alpha\right),\mathbb{F}\hat{H}\left(\sigma\right)\right\rangle -\Upsilon\left(\alpha,\alpha,\sigma\right)-\psi\left(\sigma\right).
\]
But, since $\psi\left(\sigma\right)=0$ for all $\sigma\in\hat{W}$,
\begin{equation}
\left\langle z\left(\alpha\right),\mathbb{F}\hat{H}\left(\sigma\right)\right\rangle =\Upsilon\left(\alpha,\alpha,\sigma\right),\ \ \ \forall\alpha\in T^{\ast}Q,\ \sigma\in\hat{W}.\label{sp}
\end{equation}
Thus, we have for $z$ the equations {[}see $\left(\ref{fp}\right)$
and $\left(\ref{sp}\right)${]}
\[
\left\langle z\left(\alpha\right),\mathbb{F}\hat{H}\left(\alpha\right)\right\rangle =-\mu\left(\alpha\right)\ \ \ \ \text{and}\ \ \ \left\langle z\left(\alpha\right),\mathbb{F}\hat{H}\left(\sigma\right)\right\rangle =\Upsilon\left(\alpha,\alpha,\sigma\right),
\]
or equivalently, in terms of the (non-degenerate) tensor $\hat{\mathfrak{b}}$,
we have the first two equations in $\left(\ref{fsp}\right)$. Since
$\hat{H}$ is positive-definite w.r.t. $\alpha^{\bullet}$, and consequently
$d\hat{H}\left(\alpha^{\bullet}\right)=0$, it is easy to show from
$\left(\ref{y1}\right)$ that the condition $y\left(\alpha^{\bullet}\right)=0$
implies $z\left(\alpha^{\bullet}\right)=0$. This gives us the last
part of $\left(\ref{fsp}\right)$.

Now, writing $z\left(\alpha\right)=z_{\shortparallel}\left(\alpha\right)-z_{\perp}\left(\alpha\right)$
as explained above, Eq. $\left(\ref{fsp}\right)$ translates to
\begin{equation}
\hat{\mathfrak{b}}\left(z_{\shortparallel}\left(\alpha\right),\alpha\right)=\hat{\mathfrak{b}}\left(z_{\perp}\left(\alpha\right),\alpha\right)-\mu\left(\alpha\right)\ \ \ \text{and\ \ \ }\hat{\mathfrak{b}}\left(z_{\perp}\left(\alpha\right),\sigma\right)=-\Upsilon\left(\alpha,\alpha,\sigma\right),\label{fsp1}
\end{equation}
for all $\alpha\in T^{\ast}Q$ and $\sigma\in\hat{W}$. The second
part of $\left(\ref{fsp1}\right)$ says precisely that $z_{\perp}$
is given by $\left(\ref{zo}\right)$. In particular, $z_{\perp}\left(0\right)=0$
on any fiber of $T^{\ast}Q$, so $z_{\bot}\left(\alpha^{\bullet}\right)=0$,
and consequently $z_{\shortparallel}\left(\alpha^{\bullet}\right)=0$
also. We conclude then that $z_{\shortparallel}$ must satisfies all
the conditions appearing in $\left(\ref{zos}\right)$. 

Reciprocally, it is easy to see that, if $z$ (resp. $z_{\perp}$
and $z_{\shortparallel}$) satisfies $\left(\ref{fsp}\right)$ {[}resp.
the Equations $\left(\ref{zo}\right)$ and $\left(\ref{zos}\right)${]},
and $\hat{H}$ is a solution of $\left(\ref{t0}\right)$, reversing
the steps above we have that $y$ is a solution of $\left(\ref{vincf}\right)$.\ \ \ $\square$

\bigskip{}

From all that, we have another stabilization method (included in the
LCB method).

\begin{definition} \label{slcbm}Given an underactuated system $\left(H,\mathcal{W}\right)$,
with $H$ simple and $\mathcal{W}$ defined by a subbundle $W\subset T^{\ast}Q$,
and a critical point $\alpha^{\bullet}\in T^{\ast}Q$ of $X_{H}$,
the \textbf{simple LCB method} consists in:
\begin{itemize}
\item finding a simple function $\hat{H}$ that solves $\left(\ref{t0}\right)$,
with $\hat{H}$ positive definite w.r.t. $\alpha^{\bullet}$;
\item defining $z_{\perp}$ by the Eq. $\left(\ref{zo}\right)$;
\item fixing a fiber preserving map $z_{\shortparallel}:T^{\ast}Q\rightarrow T^{\ast}Q$
by following the steps $ii$ and $iii$ of Proposition \ref{stepsiii},
replacing $\mu\left(\alpha\right)$ by $\mu\left(\alpha\right)-\hat{\mathfrak{b}}\left(z_{\perp}\left(\alpha\right),\alpha\right)$
in formula $\left(\ref{cif}\right)$;
\item defining $z:=z_{\shortparallel}-z_{\perp}$, $y$ as in Eq. $\left(\ref{y1}\right)$
and finally $Y\in\mathfrak{X}\left(T^{\ast}Q\right)$ as the vertical
lift of $y$. 
\end{itemize}
\end{definition}

\subsection{Maximality and equivalence}

The following theorem is a direct consequence of Theorem \ref{inc}
and the calculations of the two previous subsections.

\begin{theorem} \label{inc2}Consider an underactuated system $\left(H,\mathcal{W}\right)$,
with $H$ simple, $\mathcal{W}$ given by a subbundle $W\subset T^{\ast}Q$,
and a critical point $\alpha^{\bullet}\in T^{\ast}Q$ of $X_{H}$.
If we are given a vector field $Y\subset\mathcal{W}$ and a simple
Lyapunov function for $X_{H}+Y$ and $\alpha^{\bullet}$, then $Y$
is given by the simple LCB method (see Definition \ref{slcbm}).\end{theorem}

This result says that the simple LCB method is maximal among the Lyapunov
based stabilization methods for which the related Lyapunov functions
can be chosen simple. In particular, it says that the simple CH method
(see Definition \ref{ches}) is included in the simple LCB method
(recall Definition \ref{ied}). We show below that the other inclusion,
and consequently the equivalence, also holds. This is a very remarkable
fact, since the form of the control laws given by the (simple) CH
method seems to be not too general, mainly because of the rather special
form of the gyroscopic forces. But, as we show below, this is a wrong
impression.

\begin{theorem} \label{LCBeqCH} The simple LCB and CH methods are
equivalent, in the sense of Definition \ref{ied}. \end{theorem}

\emph{Proof. }Consider the set $\mathfrak{U}$ formed out by the triples
$\left(H,\mathcal{W},\alpha^{\bullet}\right)$, where $H$ is simple,
$\mathcal{W}$ is given by a subbundle $W\subset T^{\ast}Q$, and
$\alpha^{\bullet}\in T^{\ast}Q$ is a critical point of $X_{H}$.
Note that the simple LCB and CH methods are stabilization methods
on $\mathfrak{U}$ (see Definition \ref{smd}). Denote by $\digamma_{LCB}$
and $\digamma_{CH}$ their corresponding functions. Theorem \ref{inc2}
tells us that 
\[
\digamma_{CH}\left(H,\mathcal{W},\alpha^{\bullet}\right)\subset\digamma_{LCB}\left(H,\mathcal{W},\alpha^{\bullet}\right),\ \ \ \forall\left(H,\mathcal{W},\alpha^{\bullet}\right)\in\mathfrak{U}.
\]
Let us show the other inclusions also hold. Let $Y\in\digamma_{LCB}\left(H,\mathcal{W},\alpha^{\bullet}\right)$.
Then, $Y$ is the vertical lift of a fiber preserving map $y$ given
by $\left(\ref{y1}\right)$, with $\hat{H}$ simple, positive definite
w.r.t. $\alpha^{\bullet}$, solving $\left(\ref{t0}\right)$, and
$z$ satisfying $\left(\ref{fsp}\right)$. We want to see that $Y\in\digamma_{CH}\left(H,\mathcal{W},\alpha^{\bullet}\right)$.
This would mean, comparing $\left(\ref{y}\right)$ and $\left(\ref{y1}\right)$,
that 
\[
z\left(\alpha\right)=z_{d}\left(\alpha\right)+\mathfrak{Z}_{g}\left(\alpha,\alpha,\mathbb{F}\hat{H}^{-1}\left(\cdot\right)\right)
\]
for some fiber preserving map $z_{d}$ satisfying $\left(\ref{szd}\right)$
and some tensor $\mathfrak{Z}_{g}$ given by $\left(\ref{Zg}\right)$.
So, it is enough to take
\[
\left\{ \begin{array}{l}
\mathfrak{Z}_{g}\left(\alpha_{1},\alpha_{2},\sigma\right):=\Upsilon\left(\alpha_{1},\alpha_{2},\sigma\right),\ \ \ \mathfrak{Z}_{g}\left(\sigma_{1},\sigma_{2},\gamma\right):=-\Upsilon\left(\gamma,\sigma_{2},\sigma_{1}\right)-\Upsilon\left(\gamma,\sigma_{1},\sigma_{2}\right),\\
\mathfrak{Z}_{g}\left(\gamma_{1},\sigma,\gamma_{2}\right):=\mathfrak{Z}_{g}\left(\sigma,\gamma_{1},\gamma_{2}\right):=-\frac{1}{2}\Upsilon\left(\gamma_{1},\gamma_{2},\sigma\right),\ \ \ \mathfrak{Z}_{g}\left(\gamma_{1},\gamma_{2},\gamma_{3}\right):=0,
\end{array}\right.
\]
with $\alpha_{i}\in T^{\ast}Q$, $\sigma,\sigma_{i}\in\hat{W}$ and
$\gamma,\gamma_{i}\in W$, and
\[
z_{d}\left(\alpha\right)\coloneqq z\left(\alpha\right)-\mathfrak{Z}_{g}\left(\alpha,\alpha,\mathbb{F}\hat{H}^{-1}\left(\cdot\right)\right).
\]
In fact, using $\left(\ref{fsp}\right)$ and the definition of $\mathfrak{Z}_{g}$,
\[
\hat{\mathfrak{b}}\left(z_{d}\left(\alpha\right),\alpha\right)\coloneqq\hat{\mathfrak{b}}\left(z\left(\alpha\right),\alpha\right)-\mathfrak{Z}_{g}\left(\alpha,\alpha,\alpha\right)=-\mu\left(\alpha\right)+0=-\mu\left(\alpha\right),
\]
for all $\alpha\in T^{*}Q$,
\[
\hat{\mathfrak{b}}\left(z_{d}\left(\alpha\right),\sigma\right)\coloneqq\hat{\mathfrak{b}}\left(z\left(\alpha\right),\sigma\right)-\mathfrak{Z}_{g}\left(\alpha,\alpha,\sigma\right)=\Upsilon\left(\alpha,\alpha,\sigma\right)-\Upsilon\left(\alpha,\alpha,\sigma\right)=0,
\]
for all $\sigma\in\hat{W}$, and
\[
z_{d}\left(\alpha^{\bullet}\right)\coloneqq z\left(\alpha^{\bullet}\right)-\mathfrak{Z}_{g}\left(\alpha^{\bullet},\alpha^{\bullet},\mathbb{F}\hat{H}^{-1}\left(\cdot\right)\right)=0-0=0,
\]
what implies that $z_{d}$ satisfies $\left(\ref{szd}\right)$.\ \ \ $\square$

\bigskip{}

Let us mention a straightforward, but important, implication of the
Theorems \ref{inc2} and \ref{LCBeqCH} about the energy shaping method:
given a triple $\left(H,\mathcal{W},\alpha^{\bullet}\right)\in\mathfrak{U}$,
any control law $Y\subset\mathcal{W}$ such that $X_{H}+Y$ has a
simple Lyapunov function, related to the point $\alpha^{\bullet}$,
can be constructed with the simple CH method, and consequently (by
inclusion), with the (general) CH method. Roughly speaking, in the
realm of simple Hamiltonian systems with actuators, if we want to
stabilize one of them by means of finding a simple Lyapunov function,
the energy shaping method is the most general way to do it. This claim,
to the best of our knowledge, is not previously mentioned in the literature.

\bigskip{}

As a last comment, we want to say that, in spite of the equivalence,
it is not convenient to discard one of these methods. Although they
give rise to the same sets of control laws, the latter are constructed
in rather different ways. Also, the ideas behind both methods are
completely different, so, one of them can be more appropriate than
the other in certain situations. For instance, in Ref. \cite{gsz}
we were able to improve a result of \cite{chang2} about stabilizability
of underactuated Hamiltonian systems with two degrees of freedom.
More precisely, in Reference \cite{chang2}, by using the energy shaping
method, a set of conditions that ensures the stabilizability of such
systems has been established. In Reference \cite{gsz}, by using the
LCB method, we have shown that such conditions not only ensure the
stabilizability, but also the asymptotic stabilizability. 

\section*{Acknowledgements}

S. Grillo and L. Salomone thank CONICET for its financial support.

\end{document}